\documentclass[a4paper,preprint,10pt,maxbibnames=20,giveninits=true,sort&compress]{elsarticle}

\usepackage[bottom=2cm,left=2cm,right=2cm,top=2cm]{geometry}
\usepackage{amssymb}
\usepackage{xcolor}
\usepackage{float}
\usepackage[utf8]{inputenc}
\usepackage[T1]{fontenc}
\usepackage[bbgreekl]{mathbbol}
\usepackage{amsmath,amsfonts,amsthm}
\usepackage{bbold}
\DeclareSymbolFontAlphabet{\mathbb}{AMSb}
\DeclareSymbolFontAlphabet{\mathbbl}{bbold}
\DeclareFontFamily{OT1}{pzc}{}
\DeclareFontShape{OT1}{pzc}{m}{it}{<-> s * [1.10] pzcmi7t}{}
\DeclareMathAlphabet{\mathpzc}{OT1}{pzc}{m}{it}
\usepackage{mathtools}
\usepackage{stackengine,scalerel}
\usepackage{accents}
\usepackage{graphicx}
\usepackage[colorlinks=false, pdfborder={0 0 0}]{hyperref}
\usepackage{csquotes}
\usepackage{subcaption}
\usepackage{accents}
\usepackage{textcomp}
\usepackage{import}
\usepackage{caption,graphicx,newfloat}
\usepackage[symbol]{footmisc}
\renewcommand{\thefootnote}{\fnsymbol{footnote}}
\usepackage{tablefootnote}
\usepackage{wrapfig}

\DeclareFloatingEnvironment[
fileext=lob,
listname={List of Boxes},
name=Box,
placement=htp,
]{BOX}

\stackMath
\newcommand\tenq[2][1]{%
	\def\useanchorwidth{T}%
	\ifnum#1>1%
	\stackunder[0pt]{\tenq[\numexpr#1-1\relax]{#2}}{\scriptscriptstyle\sim}%
	\else%
	\stackunder[1pt]{#2}{\scriptscriptstyle\sim}%
	\fi%
}
\newcommand\ProjectionTensDev{\tenq[2]{P}_\mathrm{d}}
\newcommand\ProjectionTensSph{\tenq[2]{P}_\mathrm{s}}

\journal{arXiV}
\begin{document}
	\begin{frontmatter}
        \title{Empirical Hyper Element Integration Method (EHEIM) with Unified Integration Criteria for Efficient Hyper Reduced FE\textsuperscript{2} Simulations}
		\author[Freiberg]{Nils Lange}
		\author[Freiberg]{Geralf Hütter\footnote[2]{Corresponding author. \textit{E-mail address:} \href{Geralf.Huetter@imfd.tu-freiberg.de}{Geralf.Huetter@imfd.tu-freiberg.de} }}
		\author[Freiberg]{Bjoern Kiefer}
        \affiliation[Freiberg]{organization={TU Bergakademie Freiberg, Institute of Mechanics and Fluid Dynamics},
            addressline={Lampadiusstr.~4}, 
            city={Freiberg},
            postcode={09596}, 
            country={Germany}}
		\begin{abstract}
			Numerical homogenization for mechanical multiscale modeling by means of the finite element method (FEM) is an elegant way of obtaining structure-property relations, if the behavior of the constituents of the lower scale is well understood. However, the computational costs of this so-called FE\textsuperscript{2} method are so high that reduction methods are essential. While the construction of a reduced basis for the microscopic nodal displacements using proper orthogonal decomposition (POD) has become a standard technique, the reduction of the computational effort for the projected nodal forces, the so-called hyper reduction, is an additional challenge, for which different strategies have been proposed in the literature.
            The empirical cubature method (ECM), which has been proven to be very robust, implemented the conservation of the total volume is used as a constraint in the resulting optimization problem, while energy-based criteria have been proposed in other contributions.
            
            The present contribution presents a unified integration criteria concept, involving the aforementioned criteria, among others. 
             These criteria are used both with a Gauss point-based as well as with an element-based hyper reduction scheme, the latter retaining full compatibility with the common modular finite element framework.
             The methods are combined with a previously proposed clustered training strategy and a monolithic solver. Numerical examples empirically demonstrate that the additional criteria improve the accuracy for a given number of modes. Vice verse, less modes and thus lower computational costs are required to reach a given level of accuracy.  
             
		\end{abstract}
		\begin{keyword}
			Computational Homogenization\sep Multiscale Simulation \sep FE\textsuperscript{2}\sep Reduced Order Modeling \sep Hyper Reduction\sep Unified Criteria
		\end{keyword}
	\end{frontmatter}

    
\section{Introduction}

The overall behavior of engineering structures is largely governed by the design of their materials, where the respective characteristic material heterogeneities and defects at a certain length scale interact. An accurate description of the material behavior is therefore crucial for assessing the stiffness, strength and integrity of structures under complex loading conditions. When a pronounced scale separation between the structure's and the material's length scales is present, multiscale modeling can be used to couple the scales and thereby combine the advantage of the microscopic model's accuracy with the macroscopic model's efficiency \cite{kanoute_multiscale_2009}.

For the description of the overall behavior, a broad variety of analytical and semi-analytical methods has been developed. These mean-field theories are mainly based on \textsc{Eshelby}'s equivalent inclusion problem \cite{Eshelby_1957,Mura_1987} and require only statistical descriptions of the microstructural constituents. A profound drawback of such approaches is their low flexibility and generality, especially for complex geometries and irreversible behavior, albeit incremental approximations exist for elasto-plastic \cite{Tandon_1988} and other physically and also geometrically non-linear behaviors. In order to come to more accurate descriptions, the concept of the Representative Volume Element (RVE)---a term first used by \citet{Hashin_1963}---can be introduced, which refers to a material volume containing a sufficient number of relevant micro-constituents and is, in an averaged sense, representative of the material's continuum properties on the macroscopic scale \cite{NematNasser_1993}. When the Finite Element Method (FEM) is utilized at both scales, the conventionally expensive \enquote{FE\textsuperscript{2} method} emerges, a method used for the first time in the late 1990s, see \cite{Feyel_1999,Miehe_1999,Smit1998}. In order to lower the computational expenses, a variety of numerical schemes have been developed that require fewer resources for the evaluation of the RVE response. 
A detailed review can be found in \citet{Geers_multiscale_2017}.
Nowadays, machine learning based surrogate models have become immensely popular, cf.~\citet{Bishara_2022}, which, however, conceptually lead away from the general FE\textsuperscript{2} setting, often at the expense of flexibility, accuracy and/or generalization capabilities.

In contrast, reduced-order models (ROM), which have already been very successfully  employed in various areas of engineering and computer science, retain the fundamental physics-based structure of the equations, while significantly reducing the number of unknowns.
Soon after the millennial turn, ROM has been applied to the FE\textsuperscript{2} problem \cite{Ryckelynck_2005,Yvonnet2007,Hernandez_2014,Farhat2014}, while other groups \cite{Oskay2007,Fritzen2013} incorporated ROM to the eigendeformation fields within transformation field analysis (TFA)-type approaches.
Here, we focus on the former-type, where a low dimensional representation of the nodal displacement vector is employed, due to its larger flexibility.
The applications have not been limited to mechanical \cite{Rocha_2019,Bhattacharjee2020,Raschi_2021,Wulfinghoff_2024,Scheunemann2024,Guo2024,Lange_hyper_ROM_2024,Soldner2017} and multi-field \cite{Brands2019,Zhang2025,Wulfinghoff2025a,Moreau2025} computational homogenization, but also applications to full-scale simulations have been reported \cite{Kerfriden2011,Hernandez_2016,Oezmen2021}, both mostly using Galerkin projection.
However, computing FEM-type integral expressions necessary for solving the weak form (i.e., computing the internal force $\underline{\hat{f}}_\mathrm{int}$) is still expensive. For approximating these integrals, so-called hyperreduction approaches, a term coined by \citet{Ryckelynck_2005}, were developed. They can be classified into interpolation methods and cubature methods \cite{Hernandez_2016,Hernandez_2020,Tuijl_2018}. Comparative studies of Mergheim, Steinmann and co-workers \cite{Soldner2017,Brands2019} have shown that cubature methods are the more robust schemes and preserve fundamental properties such as symmetry and positive definiteness of the stiffness matrix, if present, which is the reason we focus on this kind of methods.
Their main idea is to find an underintegration scheme that preserves certain field variable averages for known solutions gathered in an offline phase, with the intention that in the online simulation phase, these averages are still sufficiently close to the
ones that would have been computed with full integration. The Empirical Cubature Method (ECM) of \citet{Hernandez_2016} is founded on the direct integration of the ROM-reduced internal force $\underline{\tilde{f}}_\mathrm{int}$, as is the offline-free scheme presented by \citet{Rocha_2019}, while other works are founded on the integration of the elastic free energy $\langle\,\psi^\mathrm{el}\,\rangle_{\mathcal{V}}$ \cite{Oliver2017,Caicedo_2019,Raschi_2021}, or on the strain $\langle\,\tenq{\varepsilon}\,\rangle_{\mathcal{V}}$ and strain fluctuation covariance $\langle\,\left[\tenq{\varepsilon}-\langle\,\tenq{\varepsilon}\,\rangle_{\mathcal{V}}\right]\otimes\left[\tenq{\varepsilon}-\langle\,\tenq{\varepsilon}\,\rangle_{\mathcal{V}}\right]\,\rangle_{\mathcal{V}}$ \cite{Wulfinghoff_2024}. In this article, a unified concept shall be introduced that allows for several of these criteria to be included, and its performance is tested in representative numerical examples. Another aspect is that the idea of cubature methods is not fully compatible with the modular structure of the FEM, since access of the solver to particular hyper integration points inside elements is required.
Hence, an approach is used in the present contribution which does not rely on finding hyper integration points, but instead builds on the concept of \enquote{hyper elements} that preserves compatibility with the standard finite element framework, cf.~\cite{Farhat2014,Brands2019}.
\section{Theory}
\subsection{Mechanical two-scale boundary value problem}
Let us consider the classical \enquote{first-order} mechanical, quasi-static, two-scale boundary value problem. For convenience, the equations are first provided in the small deformation setting, before the necessary modifications for the large deformation regime are briefly mentioned. Mechanical equilibrium is then expressed as
\begin{equation}\label{eqn:BLM}
	\underline{\nabla}_{\;\!\underline{X}}\cdot\tenq{\Sigma}+\underline{B}=\underline{0}\ \ \ \text{and}\ \ \ \underline{\nabla}_{\;\!\underline{x}}\cdot\tenq{\sigma}=\underline{0}\ ,\ \ 
\end{equation}
where here, and in the following, corresponding variables on the micro- and macro-scale are written with the same symbol, but in lower and upper case, respectively. On both scales, displacement fields exist, $\underline{U}(\underline{X})\in\mathbb{R}^d$ and $\underline{u}(\underline{x})\in\mathbb{R}^d$ ($d\in\{2,3\}$), through which the corresponding strain measures are defined as
\begin{equation}\label{eqn:strains}
\tenq{E}:=\mathrm{sym}(\mathrm{grad}(\underline{U}))\ ,\ \ \ \tenq{\varepsilon}:=\mathrm{sym}(\mathrm{grad}(\underline{u})) \ .
\end{equation}
The Hill-Mandel lemma \cite{Hill_1967}, also known as condition of macro-homogeneity, states that the stress power produced at the macroscale $P$ must be equal to the averaged microscopic counterpart $p$
\begin{equation}
    P:=\tenq{\Sigma}:\tenq{\dot{E}}\stackrel{!}{=}\dfrac{1}{V}\int_\mathcal{V}\tenq{\sigma}:\tenq{\dot{\varepsilon}}\ \mathrm{d}\;\!V=\langle p\rangle_\mathcal{V}\ ,
    \label{eq:Hill-Mandel}
\end{equation}
wherein $\mathcal{V}$ represents the RVE domain and $V=\mathrm{vol}(\mathcal{V})$ its volume. This condition can be fulfilled by defining the macroscopic stress $\tenq{\Sigma}$ as the volume-averaged microscopic stress $\tenq{\sigma}$
\begin{equation}
    \tenq{\Sigma}=\dfrac{1}{V}\int_\mathcal{V}\tenq{\sigma}(\tenq{\varepsilon},\mathbbl{a})\ \mathrm{d}\;\!V=\langle\tenq{\sigma}\rangle_\mathcal{V}\ .
    \label{eq:Sigmamacro}
\end{equation}
Here, the micro-stress state is a function of the current strain state $\tenq{\varepsilon}$ and a set of history variables $\mathbbl{a}$. This assumes  assuring the equality of the macroscopic strain $\tenq{E}$, as introduced in \eqref{eqn:strains}, and the respective micro-average
\begin{equation}\label{eqn:strain_equiv}
    \tenq{E}=\dfrac{1}{V}\int_\mathcal{V}\tenq{\varepsilon}\ \mathrm{d}\;\!V\ .
\end{equation}
In this work, periodic boundary conditions are applied in order to fulfill condition \eqref{eqn:strain_equiv}, which results in the linear constraint
\begin{equation}\label{eqn:PBC}
    \underline{u}^+=\underline{u}^-+\tenq{E}\cdot\left[\underline{x}^+-\underline{x}^-\right]
\end{equation}
for homologeous points $\underline{x}^+$ and $\underline{x}^-$ of the RVE boundary. Details on how to derive this equation can, for instance, be found in the article by \citet{Miehe_2002}.

In the large deformation regime, the first Piola-Kirchhoff $\tenq{\Sigma}^\mathrm{PK1}$ stress and the deformation gradient $\tenq{F}=\tenq{R}\cdot\tenq{U}$ are instead used in these equations as the respective measures of stress and deformation.
The rigid body rotations $\tenq{R}$ can be eliminated from the micro-scale, and thus from the reduced-order model, by transferring only the right stretch tensor $\tenq{U}$ to the micro-scale, i.e., substituting $\tenq{E}=\tenq{U}-\tenq{I}$ into the equations above. Then, the symmetric Biot stress $\tenq{\Sigma}^\mathrm{B}=\mathrm{sym}(\tenq{R}^\mathrm{T}\cdot\tenq{\Sigma}^\mathrm{PK1})$ arises as the work-conjugate quantity $\tenq{\Sigma}^\mathrm{B}:\dot{\tenq{U}}=\tenq{\Sigma}^\mathrm{PK1}:\dot{\tenq{F}}$, see \ \cite{Lange_hyper_ROM_2024,Kochmann_2018,Kunc2019}. 

\subsection{Finite element formulation}
Finally, both problems shall be solved with the FEM, for which the weak forms of equation \eqref{eqn:BLM} are needed. Because the main focus in this paper lies on the approximation of the microscopic problem, in the following, the proper treatment of the macroscopic problem is assumed to be  given. The weak form of the microscopic formulation reads
\begin{equation}\label{eqn:weak}
    \delta\;\!w=\int_{\mathcal{V}}\tenq{\sigma}:\delta\;\!\tenq{\varepsilon}\ \mathrm{d}\;\!V -\delta\;\!w_\mathrm{ext}=0 \ ,
\end{equation}
wherein the external (virtual) work $\delta\;\!w_\mathrm{ext}=\delta\;\!\tenq{E}:\tenq{\Sigma}\,V$ arises from the Hill-Mandel condition \eqref{eq:Hill-Mandel}.
 Through the spatial discretization of set $\mathbbl{z}$ into $m_\mathrm{ele}$ finite elements, with element domain $\mathcal{V}_e$ and $\mathcal{V}\approx\bigcup_{e=1}^{m_\mathrm{ele}}\mathcal{V}_e\nonumber$, an element-wise approximation of the displacement field $\underline{u}(\underline{x})$, for $\underline{x}\in\mathcal{V}_e$, is introduced by means of the nodal displacement coefficients $\underline{\hat{u}}_{\;\!e}$ as
\begin{equation}
    {\underline{u}(\underline{x},t)=\underline{\underline{N}}(\underline{x})\cdot\underline{\hat{u}}_{\;\!e}(t)}\ ,
\end{equation}
where  $\underline{\underline{N}}$ denotes the matrix of shape functions. Likewise, the same approximation is used for the test field $\delta\,\underline{u}=\underline{\underline{N}}\cdot\delta\underline{\hat{u}}$ (Galerkin approach), with arbitrary test field coefficients $\delta\underline{\hat{u}}$. 
The strain-displacement relation  $[\;\!\tenq{\varepsilon}\;\!]=\underline{\underline{B}}\cdot\underline{{\hat{u}}}_{\;\!e}$ can thus be written by means of the $\underline{\underline{B}}$-matrix, where $[\;\!\tenq{\bullet}\;\!]$ indicates the Voigt representation of a tensor.
Inserting these expressions into the weak form, i.e., 
\begin{equation}
    \delta\;\!w_\mathrm{int}=\sum_{e\in\mathbbl{z}}{\delta\underline{\hat{u}}_e}^\mathrm{T}\cdot\underbrace{\int_{\mathcal{V}_e}\underline{\underline{B}}^\mathrm{T}\cdot[\;\!\tenq{\sigma}\;\!]\,\mathrm{d}\;\!V}_{\textstyle:=\underline{\hat{f}}_\mathrm{\;\!int}^{\;\!e}}
\end{equation}
allows identifying the element nodal force vector $\underline{\hat{f}}_\mathrm{\;\!int}^{\;\!e}$ as indicated.
The connection of the element nodal displacements $\underline{\hat{u}}_{\;\!e}$ and the global ones $\underline{\hat{u}}$ can formally  be  written as
\begin{equation}\label{eqn:node_to_global}
    \underline{\hat{u}}_{\;\!e}=\underline{\underline{A}}_{\;\!e}\cdot\ \underline{\hat{u}}\ , 
\end{equation}
using the connectivity matrix $\underline{\underline{A}}_{\;\!e}$.
In discretized form, the periodic boundary conditions \eqref{eqn:PBC} represent a constraint on the {nodes} at the boundary, which in matrix notation is expressed as
\begin{equation}
    \underline{\hat{u}}^-=\underline{\underline{A}}_u\cdot\underline{\hat{u}}^++\underline{\underline{A}}_E\cdot[\tenq{E}]\ .
\end{equation}
Therein, $\underline{\hat{u}}^-$ comprises all dependent degrees of freedom (DOF) on that part of the boundary on which they are eliminated via the constraint, while $\underline{\hat{u}}^+$ includes all kinematically independent boundary DOFs. 
Matrix $\underline{\underline{A}}_u$ contains unit entries at the positions of the respective homologeous points, while $\underline{\underline{A}}_E$ encodes the prefactor of the relative positions $\underline{x}^+-\underline{x}^-$, from \eqref{eqn:PBC}, in matrix format.
Together with the non-boundary nodes $\underline{\hat{u}}^\mathrm{inner}$, all kinematically independent DOFs are contained in $\underline{\hat{u}}^\mathrm{ind}$. The final constraint equation that governs the relation between all nodal displacements $\underline{\hat{u}}$---and also holds for the test coefficients---, as a function of the independent DOFs $\underline{\hat{u}}^\mathrm{ind}$, and the macro-strain $\tenq{E}$ reads
\begin{equation}\label{eqn:final_constraint}
    \underline{\hat{u}}=\begin{bmatrix}
        \underline{\hat{u}}^\mathrm{inner}\\\underline{\hat{u}}^+\\\underline{\hat{u}}^-
    \end{bmatrix}
=\underbrace{\begin{bmatrix}&\\[-2.0ex]\underline{\underline{I}}&\underline{\underline{0}}\\[1.5ex]\underline{\underline{0}}&\underline{\underline{A}}_u\\[-2.5ex]&\end{bmatrix}}_{\textstyle:=\underline{\underline{A}}_u^*}\cdot\ \underline{\hat{u}}^\mathrm{ind}+\underbrace{\begin{bmatrix}\\[-2.0ex]\underline{\underline{0}}\\[1.5ex]\underline{\underline{A}}_E\\\end{bmatrix}}_{\textstyle:=\underline{\underline{A}}_E^*}\cdot[\;\!\tenq{E}\;\!]\ \ 
\ \ \ \text{with}\ \ \ 
\underline{\hat{u}}^\mathrm{ind}=\begin{bmatrix}
        \underline{\hat{u}}^\mathrm{inner}\\\underline{\hat{u}}^+
    \end{bmatrix}\ .
\end{equation}
Inserting the spatial discretization and constraint \eqref{eqn:final_constraint}, also using \eqref{eqn:node_to_global}, into weak form \eqref{eqn:weak} 
gives 
\begin{align}\label{eqn:ResidualANDStress}
    \underline{\hat{r}}(\underline{\hat{u}}^\mathrm{ind},\tenq{E}):=\sum_{e\in\mathbbl{z}}\,\underbrace{{\underline{\underline{A}}_u^*}^\mathrm{T}\cdot\underline{\underline{A}}_e^\mathrm{T}}_{\textstyle:={\underline{\underline{A}}_u^e}^\mathrm{T}}\cdot\,\underline{\hat{f}}_\mathrm{\;\!int}^{\;\!e}
    \stackrel{!}{=}\underline{0}\ ,\\ 
[\;\!\tenq{\Sigma}\;\!]=\frac{1}{V}\,\underline{f}_{\Sigma}:=\frac{1}{V}\sum_{e\in\mathbbl{z}}\,\underbrace{{\underline{\underline{A}}_E^*}^\mathrm{T}\cdot\underline{\underline{A}}_e^\mathrm{T}\cdot\,\underline{\hat{f}}_\mathrm{\;\!int}^{\;\!e}}_{\textstyle:=\underline{f}_{\Sigma}^{e}}\ ,
    \label{eqn:StressMacro}
\end{align}
from the requirement that $\delta\underline{\hat{u}}^\mathrm{ind}$ and $[\delta\tenq{E}]$ must be admissible, but are otherwise arbitrary. The second equation is the discretized version of relation \eqref{eq:Sigmamacro} for the macroscopic stress tensor $\tenq{\Sigma}$ and comprises, in essence, a sum of the nodal reaction forces $\underline{f}_{\Sigma}^{e}$ of the elements at the RVE boundary, normalized by the RVE volume $V$. Note that this definition of $[\;\!\tenq{\Sigma}\;\!]$ is work-conjugate to $[\;\!\tenq{E}\;\!]$, in the sense of the Hill-Mandel condition \eqref{eq:Hill-Mandel}, even in the FE discretization, in contrast to the widely used method of evaluating $[\;\!\tenq{\Sigma}\;\!]$ as the discretized volume average.
The spatially discretized equilibrium condition $\underline{\hat{r}}=\underline{0}$ is an implicit nonlinear equation for $\underline{\hat{u}}^\mathrm{ind}$ and has to be solved iteratively, e.g.,  using Newton's method, for which the tangent matrix
\begin{equation}\label{eqn:NewtonProced}
    \underline{\underline{k}}:=\dfrac{\mathrm{d\;\!}\underline{\hat{r}}}{\mathrm{d\;\!}\underline{\hat{u}}^\mathrm{\,ind}}=\sum_{e\in\mathbbl{z}}\,{\underline{\underline{A}}_u^e}^\mathrm{T}\cdot\underline{\underline{k}}^e\cdot\underline{\underline{A}}_u^e\ , \quad \mathrm{with}\ \ \underline{\underline{k}}^e:=\dfrac{\mathrm{d\;\!}\underline{\hat{f}}_\mathrm{\;\!int}^{\;\!e}}{\mathrm{d\;\!}\underline{\hat{u}}_\mathrm{\;\!e}}\ \ 
\end{equation}
is required. Equations \eqref{eqn:ResidualANDStress} and \eqref{eqn:NewtonProced} express that the system of equations is to mainly be built up as usual, with the only difference that the element contributions $\underline{\hat{f}}_\mathrm{\;\!int}^{\;\!e}$ and $\underline{\underline{k}}^e$ of the eliminated DOFs $\underline{\hat{u}}^-$ are to be shifted to the driving DOFs $\underline{\hat{u}}^+$ and $ [\;\!\tenq{E}\;\!]$.
\subsection{Reduced-order modeling (ROM) by proper orthogonal decomposition (POD)}
The iterative solution of the micro-scale FE problem \eqref{eqn:ResidualANDStress} is rather expensive, which becomes a problem when it is needed within the FE\textsuperscript{2} scheme. Since potential solutions to parameterized problems (here by $\tenq{E}(t)$) come hand in hand with a high amount of repetitiveness---i.e., the actual solution vectors reside in a lower-dimension solution manifold---a certain solution $\underline{u}(\underline{x})$ can be well represented in terms of a low-dimensional solution vector $\underline{\tilde{u}}$, by using certain \enquote{modes} $\underline{\phi}_{\,i}$, specifically 
\begin{equation}\label{eqn:ROM_approx}
    \underline{u}(\underline{x},t)=\sum_{i=1}^{\tilde{n}}\,\underline{\phi}_{\,i}(\underline{x})\ \tilde{u}_i(t)\ .
\end{equation}
\renewcommand{\thefootnote}{\fnsymbol{footnote}}
%
The modes are obtained from \enquote{snapshots} $\underline{\hat{u}}^\mathrm{ind}_{\,i}$ of solutions to RVE simulations for $\underline{\hat{r}}=\underline{0}$, under certain  strain paths $\tenq{E}(t)$ prescribed onto the RVE\footnote[2]{Also possible: Stress paths $\tenq{\Sigma}(t)$, with iteration for $\tenq{E}(t)$. These paths can be found, e.g., by a clustering strategy \cite{Lange_hyper_ROM_2024}.}. The snapshots are stored in a training matrix
\begin{equation}\label{eqn:DisplSnaps}
    \underline{\underline{\hat{x}}}_u=[\,\underline{\hat{u}}^\mathrm{ind}_{\,1},\,\underline{\hat{u}}^\mathrm{ind}_{\,2},\,...\,
    \underline{\hat{u}}^\mathrm{ind}_{\,p}\,]\ ,\ \ \underline{\underline{\hat{x}}}_u\in\mathbb{R}^{n\times p} \ .
\end{equation}
The idea of proper orthogonal decomposition (POD), see  \cite{Pearson_1901}, is to perform a truncated spectral decomposition of the covariance part of the training matrix, taking into account only a given number of modes with highest impact on the variance of the training data.
The procedure can be implemented by a singular value decomposition (SVD)
\begin{equation}
    \underline{\underline{\hat{x}}}_u\stackrel{\mathrm{SVD}}{\approx}\underline{\underline{\hat{\phi}}}\cdot\underline{\underline{\sigma}}_{\;\!u}\cdot\underline{\underline{v}}_{\;\!u}^\mathrm{T}\ ,\ \ \underline{\underline{\hat{\phi}}}\in\mathbb{R}^{n\times \tilde{n}}\ ,
\end{equation}
which simply means that a certain snapshot can be represented by a linear combination of $\tilde{n}$ modes $\underline{\hat{\phi}}_{\;\!i}$, contained as columns in the reduced basis $\underline{\underline{\hat{\phi}}}$. Then, nodal displacements $\underline{\hat{u}}^\mathrm{ind}$ for actual solutions are expressed as
\begin{equation}
    \underline{\hat{u}}^\mathrm{ind}=\underline{\underline{\hat{\phi}}}\cdot\underline{\tilde{u}}\ ,
    \label{eq:reducedprojection}
\end{equation}
in terms of the modal amplitudes $\underline{\tilde{u}}$.
Furthermore, complying with  \eqref{eqn:final_constraint} yields
\begin{equation}
    \underline{\hat{u}}=\underbrace{\underline{\underline{A}}_u^*\cdot\underline{\underline{\hat{\phi}}}}_{\textstyle:=\underline{\underline{\tilde{A}}}_u^*}\cdot\ \underline{\tilde{u}}+\underline{\underline{A}}_E^*\cdot[\;\!\tenq{E}\;\!]\ .\ \ 
\end{equation}
This final constraint equation, inserted into  weak form \eqref{eqn:weak}, results in  a reduced residual $\underline{\tilde{r}}$. Together with the corresponding necessary linearization, this reads
\vspace{-2ex}
\begin{align}
\underline{\tilde{r}}:=&\sum_{e\in\mathbbl{z}}\,\overbrace{\underbrace{{\underline{\underline{\tilde{A}}}_u^*}^\mathrm{T}\cdot\underline{\underline{A}}_e^\mathrm{T}}_{\textstyle:={\underline{\underline{\tilde{A}}}_u^e}^\mathrm{T}}\cdot\,\underline{\hat{f}}_\mathrm{int}^{\;\!e}}^{\textstyle:=\underline{\tilde{f}}_\mathrm{int}^{\;\!e}}=\sum_{e\in\mathbbl{z}}\,\underline{\tilde{f}}_\mathrm{int}^{\;\!e}=\underline{0}\ ,& \underline{\underline{\tilde{k}}}&:=\dfrac{\mathrm{d\;\!}\underline{\tilde{r}}}{\mathrm{d\;\!}\underline{\tilde{u}}}=\sum_{e\in\mathbbl{z}}\,{\underline{\underline{\tilde{A}}}_u^e}^\mathrm{T}\cdot\underline{\underline{k}}^e\cdot\underline{\underline{\tilde{A}}}_u^e\ .\ \ 
\end{align}
Note that each element $e$ contributes with $\underline{\tilde{f}}_\mathrm{int}^{\;\!e}$ to the residual $\underline{\tilde{r}}$, in contrast to the unrestricted FE formulation, where each element only contributes to the residual entries corresponding to the degrees of freedom of the nodes it is connected to. Further note that the computation of the macroscopic stress $\tenq{\Sigma}$ according to equation \eqref{eqn:StressMacro} still holds.

Classically, the system of equations $\underline{\tilde{r}}(\underline{\tilde{u}},\tenq{E})=0$ is solved at each macroscopic Gauss point by the Newton-Raphson method for a prescribed estimate $\tenq{E}=\underline{\underline{B}}\cdot\underline{U}$ of the macroscopic strain within a macroscopic Newton loop for $\underline{U}$. Alternatively, the system of equations for the microscopic (modal) displacements $\underline{\tilde{u}}$ and macroscopic displacements $\underline{U}$ can be solved in a \emph{monolithic} way, in a common Newton-Raphson loop (using static condensation). Details on the monolithic solution strategy, which is also used here to avoid the costly iterations at the micro-scale, are given in \cite{Gruttmann_2013,Okada_2010,Lange_2021,Lange_hyper_ROM_2024,Eidel2019}.

\subsection{Empirical Hyper Element Integration Method (EHEIM)\label{sec:EHEIM}}
The key idea of the Empirical Cubature Method (ECM), cf.~\cite{Hernandez_2016, Hernandez_2020}, is to find an underintegration scheme with much fewer integration points $\tilde{m}$ than the original scheme, i.e., $m\gg\tilde{m}$, to compute the reduced residual $\underline{\tilde{r}}$ and macroscopic stress $\tenq{\Sigma}$. The hyper integration scheme is \enquote{calibrated} such that these measures are optimally preserved for given snapshots, for whose generation the same training paths $\tenq{E}(t)$ (or $\tenq{\Sigma}(t)$) are chosen here, which were used to find the displacement snapshots, see \eqref{eqn:DisplSnaps}. 
The integration of the ECM approach into established FE frameworks is generally a problem from the implementation standpoint, as  already stated in the introduction, because in the actual FE\textsuperscript{2} simulation, the microscopic solver must act on the element level, to enforce the integral evaluations at the determined hyper integration points with the determined weights, which is not possible in the usual program structure of the FEM. To overcome this, a promising approach is to select whole \enquote{hyper elements} \cite{Farhat2014,Brands2019}, whose output is weighted with scalar factors $\gamma_e$---an approach we call \enquote{Empirical Hyper Element Integration Method} (EHEIM). This is favorable, because then the elements do not need to be modified, but standard elements can be used  to compute the projected nodal forces. Including the nodal reaction forces $\underline{f}_\Sigma$ required for computing the macroscopic stress $\tenq{\Sigma}$, from \eqref{eqn:StressMacro}, in a single step, this approach can be expressed as
\begin{align}
\underline{\tilde{f}}:=\begin{pmatrix}\underline{\tilde{r}} \\ \underline{f}_\Sigma\end{pmatrix} \approx\sum_{e\in\mathbbl{\tilde{z}}}\gamma_e\,\underline{\tilde{f}}_\mathrm{int}^{\;\!e}
=\sum_{e\in\mathbbl{\tilde{z}}}\gamma_e\,{\underline{\underline{\tilde{A}}}_u^e}^\mathrm{T}\cdot\underline{f}_\mathrm{int}^{\;\!e} \ ,
\end{align}
so that both ROM projection as well as the hyper integration can be performed at the global level of the solver.
Therein, set $\mathbbl{\tilde{z}}$ that contains $\tilde{m}_\mathrm{ele}$ hyper elements shall naturally be a subset of the original set of elements, i.e.,  $\mathbbl{\tilde{z}}\subset\mathbbl{z}$,  with $m_\mathrm{ele}\gg\tilde{m}_\mathrm{ele}$. 
This set $\mathbbl{\tilde{z}}$ of hyper elements and their weights $\gamma_e$ could be determined from minimizing the error in $\underline{\tilde{r}}$ from snapshots, under the condition that the total volume is preserved, as is done in the ECM.
However, preserving additional quantities in an averaged sense, such as the free energy $\Psi=\langle\psi\rangle$, strain $\tenq{E}=\langle\tenq{\varepsilon}\rangle$, stress power $P=\langle p\rangle$, among others, can help in finding an underintegration scheme which generalizes better to  unseen loading conditions, as will be empirically shown in Section \ref{sec:NumExp}.

To pursue this approach, let us consider a certain field variable $\omega(\underline{x})$, whose average $\Omega$ shall be preserved by the hyper integration scheme
\begin{equation}
    \Omega:=\frac{1}{V}\int_{\mathcal{V}}\omega(\underline{x})\ \mathrm{d\;\!}V\overset{\text{FEM}}{\approx}\dfrac{1}{V}\sum_{e\in\mathbbl{z}}\bar{\omega}_e\overset{\text{EHEIM}}{\approx}\dfrac{1}{V}\sum_{e\in\mathbbl{\tilde{z}}}\gamma_e\ \bar{\omega}_e\ ,\quad \text{with}\ \ \ \bar{\omega}_e:=\int_{\mathcal{V}_e}\omega(\underline{x})\ \mathrm{d\;\!}V\ .
\end{equation}
 Analogously to the ECM, the integration consistency should be respected, meaning that the volume is integrated exactly, which is sensible for any integration scheme and reads here
\begin{equation}
    V^*:=\sum_{e\in\mathbbl{z}}V_e\stackrel{!}{=}\sum_{e\in\mathbbl{\tilde{z}}}\gamma_eV_e
\end{equation}
where $V_e:=\mathrm{vol}(\mathcal{V}_e)$ is the volume of the $e$-th element. Note that the total bulk element volume $V^*$ does not coincide with the RVE volume $V$ if pores are present. For this case, the relative bulk density $V_\mathrm{rel}:=V^*/{V}$ is introduced.
The hyper elements $\tilde{\mathbbl{z}}$ and the corresponding weighting factors, which shall be stored in vector form as $\underline{\gamma}$, can be found analogously to the ECM method, through an error minimization problem
\begin{align}
  (\mathbbl{\tilde{z}},\underline{\gamma})=\arg&\min\!\bigg(\,\sum_{k=1}^{p}\bigg(\bigg\Vert\,\underline{\underline{\eta}}_{f}\cdot\sum_{e\in\mathbbl{\tilde{z}}}\!\gamma_e\,\underline{\tilde{f}}_\mathrm{int}^{\;\!e}\bigg\Vert^2_k+\bigg\Vert\,\underline{\underline{\eta}}_{\Sigma}\cdot\bigg[\sum_{e\in\mathbbl{\tilde{z}}}\gamma_e\,\underline{f}_\Sigma^{\;\!e}-V[\;\!\tenq{\Sigma}\,]\bigg]\bigg\Vert^2_k+\bigg\Vert\,\eta_p\bigg[\sum_{e\in\mathbbl{\tilde{z}}}\gamma_e\,\bar{p}_{\;\!e}-V\!\:P\,\bigg]\bigg\Vert^2_k\\
    \nonumber ...\,+&\bigg\Vert\,\underline{\underline{\eta}}_{\varepsilon}\cdot\bigg[\sum_{e\in\mathbbl{\tilde{z}}}\gamma_e\,\bar{\tenq{\varepsilon}}_{\;\!e}-V\!\;\tenq{E}\,\bigg]\bigg\Vert^2_k+\bigg\Vert\,\eta_\psi\bigg[\sum_{e\in\mathbbl{\tilde{z}}}\!\gamma_e\,\bar{\psi}_{\;\!e}-V\,\Psi\,\bigg]\bigg\Vert^2\bigg):V^*=\sum_{e\in\mathbbl{\tilde{z}}}\gamma_eV_e\land\underline{\gamma}\in\mathbb{R}_+^{\tilde{m}}\land\tilde{\mathbbl{z}}\subset\mathbbl{z}\bigg)\ .
\end{align}
Essentially, we seek to find a certain given number $\tilde{m}_\mathrm{ele}$ of hyper elements, with corresponding positive scaling factors $\gamma_e$, under the condition of an exact volume integration. The error is composed of a sum over all snapshots, with contributions from all chosen measures, e.g.,  $P=\langle p\rangle$, $\tenq{E}=\langle\tenq{\varepsilon}\rangle$ etc., which feed into the squared error between the respective hyper-integrated and fully-integrated results (upper case symbols). The contributions---since they have different physical units---are nondimensionalized by the weighing matrices $\underline{\underline{\eta}}_\bullet$, appropriately chosen later. In order to make the problem accessible for a computer implementation, we aim to find a corresponding matrix formulation. Analogously to Hern\'andez et al.~\cite{Hernandez_2014},  this can equivalently be written with the aid of a training matrix $\underline{\underline{x}}_f$ as
\begin{equation}\label{eqn:reformulatedMIN}
    (\mathbbl{\tilde{z}},\underline{\gamma})=\arg\min\bigg(\underbrace{\left\Vert\prescript{\tilde{\mathbbl{z}}}{\phantom{}}{\underline{\underline{x}}_f}\cdot\underline{\gamma}\right\Vert^2}_{\textstyle:=\mathfrak{r}}:V^*=\underline{V}^\mathrm{T}\cdot\underline{\gamma}\land\underline{\gamma}\in\mathbb{R}^{\tilde{m}}\land\tilde{\mathbbl{z}}\subset\mathbbl{z}\bigg)\ ,
\end{equation}
where each column of $\underline{\underline{x}}_f$ corresponds to an element and each row to a certain snapshot. Furthermore, \underline{V} is the vector of element volumes
\begin{equation}
    \underline{V}=\left[V_1,\,V_2,\,\dots,\,V_{m_\mathrm{ele}}\right]^\mathrm{T}\ .
\end{equation}
The training matrix $\underline{\underline{x}}_f$ is built up from submatrices 
\begin{equation}  \underline{\underline{x}}_f^\mathrm{T}=\left[\,\underline{\underline{x}}_f^1\cdot\underline{\underline{\eta}}_f^1,\,\underline{\underline{x}}_\Sigma^1\cdot\underline{\underline{\eta}}_\Sigma^1,\,...\,,\underline{x}_\psi^1\eta_\psi^1,\,...\,,\underline{\underline{x}}_f^p\cdot\underline{\underline{\eta}}_f^p,\,\underline{\underline{x}}_\Sigma^p\cdot\underline{\underline{\eta}}_\Sigma^p,\,...\,,\underline{x}_\psi^p\eta_\psi^p\,\right]\ ,
\end{equation}
which read 
\begin{equation}\label{eqn:additionalcriteriatrainingmatrix}
\underline{\underline{x}}_{{\;\!f}}^k=\begin{bmatrix}
		\phantom{.}\\[-2.3ex]
		{\underline{{\tilde{f}}}_{\;\!\mathrm{int}}^{\;\!1}}^\mathrm{\!\!\!\!T}\\
		\vdots\\
		{\underline{{\tilde{f}}}^{\;\!m_\mathrm{ele}}_{\;\!\mathrm{int}}}^\mathrm{\!T}\\[-2.3ex]
		\phantom{.}
	\end{bmatrix}_k\!\!,\;\;\underline{\underline{x}}_{{\;\!\Sigma}}^k=\begin{bmatrix}
		\phantom{.}\\[-2.3ex]
		{\underline{{{f}}}_{\;\!\Sigma}^{\;\!1}}^\mathrm{\!\!T}\\
		\vdots\\
		{\underline{{{f}}}^{\;\!m_\mathrm{ele}}_{\;\!\Sigma}}^\mathrm{T}\\[-2.3ex]
		\phantom{.}
	\end{bmatrix}_k\!\!-\dfrac{\underline{V}\cdot\,[\tenq{\Sigma}]^\mathrm{T}_k}{V_\mathrm{rel}}\, ,\;\;\underline{x}_{\;\!p}^k\!=\!\begin{bmatrix}
	\bar{p}_1\\
	\vdots\\
	\bar{p}_{m_\mathrm{ele}}\\[-2ex]
	\phantom{.}
	\end{bmatrix}_k\!\!\!-\dfrac{\underline{V}\,P_k}{V_\mathrm{rel}}\,,\;\;...\,,\;\;\underline{x}_{\;\!\psi}^k\!=\!\begin{bmatrix}
	\bar{{\psi}}_1\\
	\vdots\\
	\bar{{\psi}}_{m_\mathrm{ele}}\\[-2ex]
	\phantom{.}
	\end{bmatrix}_k\!\!\!-\dfrac{\underline{V}\,\Psi_k}{V_\mathrm{rel}}\ .
\end{equation}
The weighting matrices $\underline{\underline{\eta}}_\mathfrak{\,\bullet}^k$ are now chosen such that each row in the training matrix becomes a unit vector
\begin{equation}\label{eqn:penalty}
	\underline{\underline{\eta}}_\mathfrak{\,\bullet}^k:=\mathrm{diag}\left(\eta^k_{1\,\bullet},\eta^k_{2\,\bullet},\,\dots\,\right)\ , \quad \text{with}\ \ \ \eta^k_{j\,\bullet}=\begin{cases}
		\Vert{\underline{x}}_{j\,\bullet}^k\Vert^{-1}&\!,\ \Vert{\underline{x}}_{j\,\bullet}^k\Vert\ge\varepsilon\\
		0&\!,\ \Vert{\underline{x}}_{j\,\bullet}^k\Vert<\varepsilon
	\end{cases}\ ,
\end{equation}
which weights each snapshot equally. This normalization of the \enquote{raw data} is common in data science and is part of so-called \enquote{feature engineering} techniques \cite{Reddy_2020}.

So far, the solution of \eqref{eqn:reformulatedMIN} would require the testing of all possible combinations of the given number of $\tilde{m}$ elements, to find the $\underline{\gamma}$ that minimizes the error $\mathfrak{r}$ and fulfills the constraints. This approach, however, is quite unfavorable, since the training matrix $\underline{\underline{x}}_f$ is very large and it is therefore too expensive to test all combinations. Instead, we follow the approach of \citet{Hernandez_2014} and start with a singular value decomposition of $\underline{\underline{x}}_f$, which reads
\begin{equation}
	\underline{\underline{x}}_f^\mathrm{T}\overset{\text{SVD}}{=}\underline{\underline{\lambda}}\cdot\underline{\underline{\sigma}}_f\cdot\underline{\underline{v}}_f^\mathrm{T}\ .
\end{equation}
Substituting this representation into the error $\mathfrak{r}$  of equation \eqref{eqn:reformulatedMIN} allows simplifying the expression to \vspace{-2ex} 
\renewcommand{\thefootnote}{\fnsymbol{footnote}}
\footnotetext[2]{The notation $R(\underline{\underline{x}})$ refers to the rank of a matrix $\underline{\underline{x}}$.}
\begin{equation}\label{eqn:SVDError}
	\mathfrak{r}=\big\Vert\,{\prescript{\tilde{\mathbbl{z}}}{\phantom{}}{\underline{\underline{x}}_f}}\cdot\underline{\gamma}\,\big\Vert^2=\Big\Vert\,{\underline{\underline{v}}_f}\cdot{\underline{\underline{\sigma}}_{f}}\phantom{}^\mathrm{T}\cdot\prescript{\tilde{\mathbbl{z}}}{\phantom{}}{{\underline{\underline{\lambda}}}}\phantom{}^\mathrm{T}\!\cdot\,\underline{\gamma}\,\Big\Vert^2=\underline{\gamma}^{\mathrm{T}}\cdot\prescript{\tilde{\mathbbl{z}}}{\phantom{}}{{\underline{\underline{\lambda}}}}\cdot{\underline{\underline{\sigma}}_{f}}\cdot\underbrace{\underline{\underline{v}}_f^\mathrm{T}\cdot\underline{\underline{v}}_f}_{\textstyle=\underline{\underline{I}}}\!\cdot\,{\underline{\underline{\sigma}}_{f}}\phantom{}^\mathrm{T}\cdot\prescript{\tilde{\mathbbl{z}}}{\phantom{}}{{\underline{\underline{\lambda}}}}\phantom{}^\mathrm{T}\!\cdot\,\underline{\gamma}=\!\!\sum_{i=1}^{\mathrm{R}(\underline{\underline{\sigma}}_f)\text{\footnotemark[2]}}\!\sigma_{i\,f}^2\left[\prescript{\tilde{\mathbbl{z}}}{\phantom{}}{{\underline{\underline{\lambda}}}}\phantom{}^\mathrm{T}\cdot\underline{\gamma}\right]^2_i \ .
\end{equation}
Since the matrix of right singular vectors $\underline{\underline{v}}_f$ is a unitary matrix, the product $\underline{\underline{v}}_f^\mathrm{T}\cdot\underline{\underline{v}}_f=\underline{\underline{I}}$ drops out. 
Expression \eqref{eqn:SVDError} then allows simplifying \eqref{eqn:reformulatedMIN} to
\begin{equation}\label{eqn:min_form_w_sing_vectors}
	(\underline{\gamma},\tilde{\mathbbl{z}})
	\!=\arg\ \min
	\Bigg(\sum_{i=1}^{\mathrm{R}(\underline{\underline{\sigma}}_f)}\sigma_{i\,f}^2\left[\prescript{\tilde{\mathbbl{z}}}{\phantom{}}{{\underline{\underline{\lambda}}}}\phantom{}^\mathrm{T}\cdot\underline{\gamma}\right]^2_i\ :\ V^*=\underline{V}^\mathrm{T}\cdot\underline{\gamma}\ \land\ \underline{\gamma}\in\mathbb{R}^{\tilde{m}_\mathrm{ele}}_+\ \land\ \tilde{\mathbbl{z}}\subset\mathbbl{z}\Bigg)\ .
\end{equation}
This formulation shows those columns $\underline{\lambda}_i$ which correspond to the largest singular values $\sigma_{i\,f}$ have the highest influence on the error $\mathfrak{r}$. Consequently, the idea is to only consider the first $\tilde{m}-1$ vectors of $\underline{\underline{\lambda}}$, which yields the truncated left singular vector matrix $\underline{\underline{\tilde{\lambda}}}\in\mathbb{R}^{m\times\tilde{m}-1}$ and transforms the minimization problem \eqref{eqn:min_form_w_sing_vectors} to a formulation of the form
\begin{equation}
	\exists\  (\underline{\gamma},\tilde{\mathbbl{z}}):\!\begin{bmatrix}
			\phantom{.}\\[-2ex]
			\prescript{\tilde{\mathbbl{z}}}{\phantom{}}{{\underline{\underline{\tilde{\lambda}}}}^\mathrm{T}}\\[1.5ex]
			\underline{V}^\mathrm{T}\\[-2ex]
			\phantom{.}
	\end{bmatrix}
	\cdot\,\underline{\gamma}=\!
	\begin{bmatrix}
			\underline{0}\\[1.5ex]
			V^*
	\end{bmatrix}
	\ \land\  \underline{\gamma}\!\in\!\mathbb{R}^{\tilde{m}_\mathrm{ele}}_+\ \land\,\tilde{\mathbbl{z}}\!\subset\!\mathbbl{z}\ .
\end{equation}
Here, the scaling factors $\underline{\gamma}$ are the solution of the shown system of equations. Choosing the elements for $\tilde{\mathbbl{z}}$ can be accomplished by following the \enquote{greedy algorithm} described in Hern\'andez et al.~\cite{Hernandez_2014}. Note that the solution is not unique,  in general, because there can be different combinations of hyper elements such that the system of equations is solvable and fulfills the constraint.

This method was implemented into the existing open-source code \texttt{MonolithFE\textsuperscript{2}}\footnote[1]{The software can be downloaded at \url{https://www.tu-freiberg.de/monolithfe2}.}, which relies on \texttt{Abaqus}~\cite{abaqus} as the macrocopic solver. Details on the implementation can be found in the preceding article \cite{Lange_hyper_ROM_2024}.

\subsection{Sampling strategy}

Both, the computation of the reduced basis as well as of the hyper integration weights and locations require snapshots from RVE simulations. The choice of loading paths $\tenq{E}(t)$ (or $\tenq{\Sigma}(t)$) for training is an important \enquote{degree of freedom} in generating a reduced-order model, both in terms of accuracy and computational costs. In the worst case, the computational costs for the respective large number of offline RVE simulations may feed off the computational gain in the online phase. 
In most studies the whole strain space $\tenq{E}$ is sampled unspecifically for the training, with strategies differing in detail \cite{Logarzo2021,Fritzen2019,Brands2019,Moreau2025}. Recently, a sampling strategy was proposed by the authors \cite{Lange_hyper_ROM_2024} in which the training paths are derived from (cheap) linear-elastic simulations of the considered macroscopic problem. The obtained paths $\tenq{E}(t)$ are clustered using the k-means algorithm and only the centroids of each cluster are used for the training. Together, they form representatives of the manifold of expected loading paths. Compared to to the unspecific sampling, this strategy does not only reduce the costs for training, but also the numbers of significant modes and thus hyper integration points can be reduced, as will be discussed also in the context of the examples below.


\section{Numerical Examples\label{sec:NumExp}}

In the following, the proposed approach is applied to three example problems to benchmark its performance and accuracy. The first one is a 2D elastic-plastic composite within small deformation theory from literature, for which the different methods are compared in detail. The second example is a porous strip, for which large deformation effects are investigated within the widely used hyper-elastic large deformation framework. Finally, a three dimensional woven composite structure from literature is considered to assess the limits regarding computational costs in industrial applications. 

\subsection{Porous elasto-plastic composite\label{sec:numerical_exp1}}

\begin{figure}[!h]
	\centering
	\subfloat[][Cantilever beam under bending with prescribed tip displacement, depth $t$ (with $t=a$).\label{Beam}]{\includegraphics[width=0.43\linewidth]{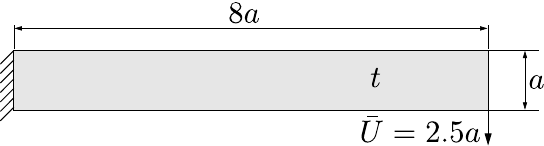}}
	\hspace{2ex}
	\subfloat[][Relative dimensions of the RVE \cite{Miehe_2002}.\label{MieheRVE}]{\includegraphics[width=0.38\linewidth]{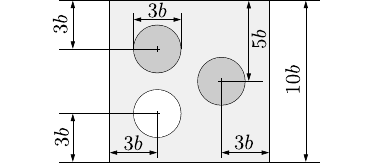}}
	\caption{Macroscopic beam problem accounting for a micro-heterogeneous material. The RVE of the microstructure comprises a soft elasto-plastic matrix, stiff elastic inclusions and a pore\label{fig:Exp1}, cf., the article by \citet{Miehe_2002}. The following relative material properties have been assumed: $E^{\mathrm{M}},\, \sigma^{\mathrm{M}}_{\mathrm{y}\,0}\!=\!0.01E^{\mathrm{M}},\,h^{\mathrm{M}}\!=\!0.016E^{\mathrm{M}}$, $E^{\mathrm{I}}\!=\!10\,E^{\mathrm{M}}$,\,$\nu^{\mathrm{M}}\!=\!\nu^{\mathrm{I}}\!=\!0.3$. Here the superscripts M and I refer to the properties of the matrix and inclusions, respectively \cite{Lange_hyper_ROM_2024}.
	}
\end{figure}

\begin{figure}[!h]
	\centering
	\subfloat[][Normalized force-displacement curves for various  $\tilde{m}$.\label{FU_HP}]{\includegraphics[width=0.43\linewidth]{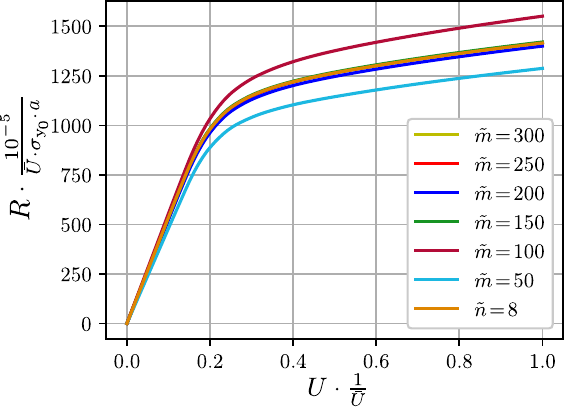}}
	\hfill
	\subfloat[][Rel.~error w.r.t.~the ROM solution for variation of $\tilde{m}$.\label{error_HP_wrt_HF}]{\includegraphics[width=0.43\linewidth]{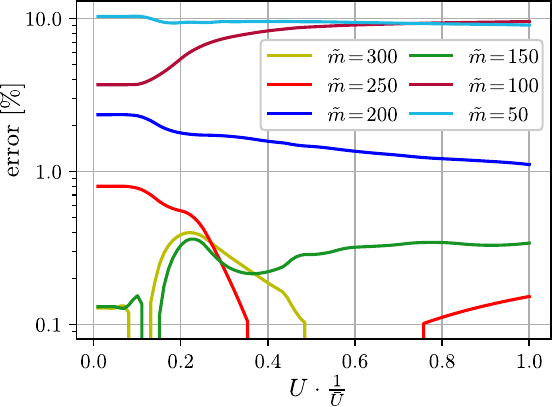}}
	\caption{Hyper integration point convergence study for a fixed number of $\tilde{n}=8$ ROM modes, using the conventional criteria ($\underline{\tilde{f}}^\mathrm{int}$ and $\underline{f}_\Sigma$). Figure taken from publication \cite{Lange_hyper_ROM_2024}.}
	\label{fig:HPconvergence}
\end{figure}

Firstly, a plane RVE of the porous elasto-plastic composite shown in Figure~\ref{fig:Exp1} is considered, as it has previously been used by \citet{Miehe_2002} for FE$^2$ investigations, and by the authors in \cite{Lange_hyper_ROM_2024}, using reduced-order modeling with ECM.
The latter investigations revealed that a (fully-integrated) ROM model, with $\tilde{n}=8$ modes, is sufficient to match the high-fidelity (HF) FE\textsuperscript{2} results with sufficient accuracy.
The results for the reaction force $R$ using ECM hyperintegration, as reported in \cite{Lange_hyper_ROM_2024}, are shown in Figure~\ref{fig:HPconvergence}. These relied on the \enquote{conventional} criteria, namely the internal force of the modes $\underline{\tilde{f}}^\mathrm{int}$ and the macroscopic stress, in form of the FE surface integral $\underline{f}_\Sigma$. The error w.r.t.\ the fully-integrated solution is observed to decrease with an increasing number of hyper integration points. However, in this case $\tilde{m}=250$ hyper integration points were necessary to reach an error below 1\%, in terms of the computed macroscopic force-displacement response at the beam tip.

\begin{figure}[!h]
\hspace{0.065\textwidth}
	\centering
	\includegraphics[width=0.43\linewidth]{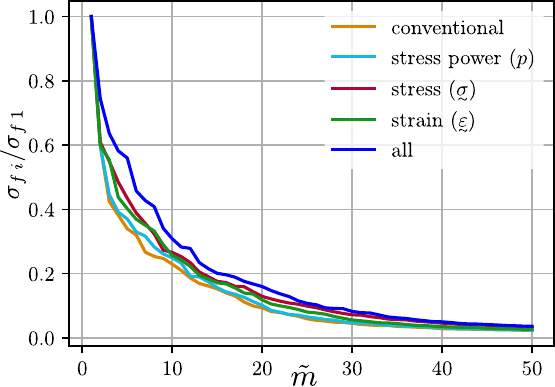}
	\caption{Decay of the singular values of the training matrix $\underline{\underline{x}}_f$ for different criteria, cf., Section \ref{sec:EHEIM}. Note that the conventional criterion was also always included in the training matrix.\label{fig:DecaySV}}
\end{figure}

This study is now repeated here, but with the additional criteria for the hyper integration selection process in the training matrix $\underline{\underline{x}}_f$, as previously described in Section \ref{sec:EHEIM}. At first, we incorporate a certain criterion, e.g.,  $\langle\,\tenq{\sigma}\,\rangle$ or $\langle\,p\,\rangle$, into $\underline{\underline{x}}_f$, while keeping the conventional criterion activated. Then, we compute the SVD and  analyze the decay of the singular values, as shown in Figure~\ref{fig:DecaySV}. 

\begin{figure}[!h]
	\centering
	\subfloat[][Normalized force-displacement curves for various  $\tilde{m}$.]{\includegraphics[width=0.43\linewidth]{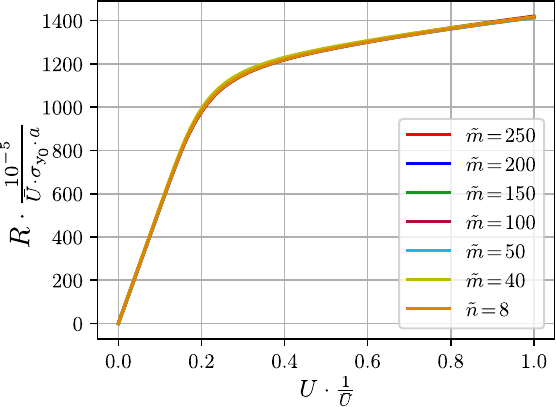}}
	\hfill
	\subfloat[][Rel.~error w.r.t.~the ROM solution for variation of $\tilde{m}$.]{\includegraphics[width=0.43\linewidth]{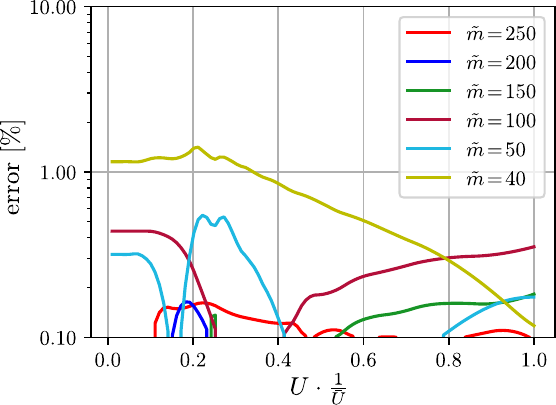}}
	\caption{Hyper integration point convergence for a fixed number of $\tilde{n}=8$ ROM modes, when enforcing the criteria of Figure~\ref{fig:DecaySV}\label{fig:IPConvergAddCrit}.}
\end{figure}

With the additional criteria, a slower decay of the singular values is observed. From this, we can conclude that the additional criteria introduce physical information into $\underline{\underline{x}}_f$ which is not simply redundant to the information already included.
The respective results for the load-deflection curve of the composite beam are shown in Figure~\ref{fig:IPConvergAddCrit}. 
Compared to the conventional ECM results in Figure~\ref{fig:HPconvergence}, it is evident that incorporating the additional criteria reduces the error by approximately an order of magnitude, for a given number $\tilde{m}$ of hyper integration points. Consequently, only $\tilde{m}=40$ to $50$ points are necessary to reach an accuracy of 1\%, in this case. This provides empirical evidence that the additional criteria are gainful in finding a hyper integration scheme that generalizes better from the training data to the actual multiscale simulation. 

\begin{figure}[ht!]
	\centering
	\includegraphics[width=0.9\linewidth]{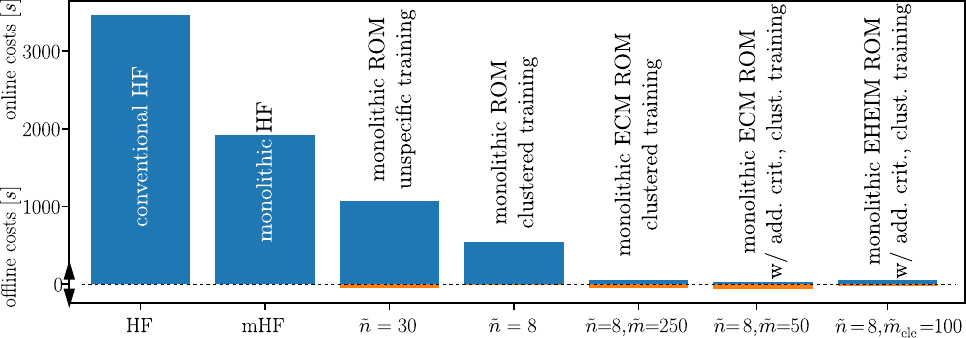}
	\caption{Computational online and offline efforts associated with the different numerical schemes.}
	\label{fig:simulation_time_example_1}
\end{figure}

The respective computational costs are compared in Figure~\ref{fig:simulation_time_example_1} and listed in detail in Table~\ref{table:simulation_times_composite}. 

\newcommand{\SpA}{0.073\linewidth}
\newcommand{\SpB}{0.06\linewidth}
\newcommand{\SpC}{0.06\linewidth}
\newcommand{\SpD}{0.06\linewidth}
\newcommand{\SpE}{0.06\linewidth}
\newcommand{\SpF}{0.06\linewidth}
\newcommand{\SpG}{0.06\linewidth}
\newcommand{\SpH}{0.06\linewidth}

\begin{table}
	\centering
	\caption{Simulation effort: The ROM offline time includes simulation of the beam with a surrogate model, RVE-FE training  simulations, and SVD on $\underline{\underline{\hat{x}}}_u$; the hyper offline time includes ROM-RVE-FE training  simulations, SVD on $\underline{\underline{x}}_f$, and hyper integration point/element selection. 
    Note that slight deviations to publication \cite{Lange_hyper_ROM_2024} in the first columns are due to the usage of different hardware and an improved version of the MonolithFE\textsuperscript{2} package, including an important efficiency-improving bugfix.
    \label{table:simulation_times_composite}}
	\begin{tabular}{|p{\SpA}|p{\SpB}|p{\SpC}|p{\SpC}|p{\SpD}|p{\SpE}|p{\SpF}|p{\SpG}|}
		\hline
		method&\multicolumn{2}{c}{HF}&\multicolumn{5}{|c|}{ROM}\\\hline
		SOE& \multicolumn{1}{c}{staggered}&\multicolumn{6}{|c|}{monolithic}\\\hline
		training& \multicolumn{2}{c}{--}&\multicolumn{1}{|c}{unspecific}&\multicolumn{4}{|c|}{clustered}\\\hline
		integr.& \multicolumn{4}{c|}{full}&\multicolumn{2}{c|}{ECM}&\multicolumn{1}{c|}{EHEIM}\\\hline
		criteria& \multicolumn{4}{c|}{--}&\multicolumn{1}{c|}{conventional}&\multicolumn{2}{c|}{additional}\\\hline
		DOFs & \multicolumn{2}{c|}{16\,385\,400\,(100\%)} & \multicolumn{1}{c|}{54\,000\,(0.3\%)}&\multicolumn{4}{c|}{14\,400\,(0.1\%)}\\\hline
		IPs & \multicolumn{4}{c|}{12\,209\,400\,(100\%)} &\multicolumn{1}{c|}{450\,000\,(4\%)}&\multicolumn{1}{c|}{90\,000\,(1\%)}&\multicolumn{1}{c|}{540\,000\,(4\%)}\\\hline
		\multicolumn{1}{|c|}{\parbox{\SpA}{\phantom{.}\vspace{-1.0ex}\\online\\time\\[-1.7ex]\phantom{.}}} & \multicolumn{1}{c|}{\parbox{\SpB}{\centering\phantom{.}\vspace{-1.0ex}\\3\,468\,s\\(100\%)}}  & \multicolumn{1}{c|}{\parbox{\SpC}{\centering\phantom{.}\vspace{-1.0ex}\\1\,915\,s\\(55\%)}}& \multicolumn{1}{c|}{\parbox{\SpD}{\centering\phantom{.}\vspace{-1.0ex}\\1070\,s\\(31\%)}}  & \multicolumn{1}{c|}{\parbox{\SpE}{\centering\phantom{.}\vspace{-1.0ex}\\542\,s\\(16\%)}} & \multicolumn{1}{c|}{\parbox{\SpF}{\centering\phantom{.}\vspace{-1.0ex}\\58\,s\\(1.7\%)}}& \multicolumn{1}{c|}{\parbox{\SpG}{\centering\phantom{.}\vspace{-1.0ex}\\26\,s\\(0.7\%)}}& \multicolumn{1}{c|}{\parbox{\SpH}{\centering\phantom{.}\vspace{-1.0ex}\\48\,s \\(1.4\%)}}\\\hline
		\multicolumn{1}{|c|}{\parbox{\SpA}{\phantom{.}\vspace{-1.0ex}\\offline\\ROM\\[-1.7ex]\phantom{.}}} & \multicolumn{2}{c|}{--}  & \multicolumn{1}{c|}{\parbox{\SpD}{\centering\phantom{.}\vspace{-1.0ex}\\47\,s\\(1.3\%)}}& \multicolumn{4}{c|}{\parbox{\SpE+\SpF+\SpE+\SpF}{\centering\phantom{.}\vspace{-1.0ex}\\13\,s\\(0.4\%)}}\\\hline
		\multicolumn{1}{|c|}{\parbox{\SpA}{\phantom{.}\vspace{-1.0ex}\\offline\\hyper\\[-1.7ex]\phantom{.}}} & \multicolumn{4}{c|}{--}  & \multicolumn{1}{c|}{\parbox{\SpD}{\centering\phantom{.}\vspace{-1.0ex}\\32\,s\\(0.9\%)}}&
		\multicolumn{1}{c|}{\parbox{\SpD}{\centering\phantom{.}\vspace{-1.0ex}\\54\,s\\(1.6\%)}} & 
		\multicolumn{1}{c|}{\parbox{\SpD}{\centering\phantom{.}\vspace{-1.0ex}\\16\,s\\(0.5\%)}}\\\hline
	\end{tabular}
    \footnotetext{Extrapolation due to a new software version for the EHEIM implementation}
\end{table}
Therein, the first bar in the plot (column in the table) corresponds to the conventional FE$^2$ with staggered Newton loops for the micro- and macro-scales, termed \enquote{high-fidelity model} (HF) in this context, and serves as a reference. 
If both scales are solved monolithically in a common loop, the computational costs reduce to 55\%.
For the reduced-order model, $\tilde{n}=30$ modes are required with an unspecific training strategy. In contrast, $\tilde{n}=8$ modes suffice to reach the same accuracy with the clustered training strategy. With this approach, the online costs drop to 16\%, if full integration is employed. Using ECM hyperintegration with the original training criteria reduces those costs by a factor of about ten, but requires some more offline time to identify the hyperintegration points and their weights.
Using the proposed additional criteria for the selection of the hyperintegration points, their number can be lowered, as discussed above, and the online costs of the multi-scale simulations are again reduced. Computing the SVD of the larger training matrix $\underline{\underline{x}}_f$ on the other hand gets slightly more expensive. This, however, might be of minor relevance for more complex macroscopic problems.

\begin{figure}[!h]
	\centering
	\subfloat[][Normalized force-displacement curves for various  $\tilde{m}_\mathrm{ele}$.]{\includegraphics[width=0.47\linewidth]{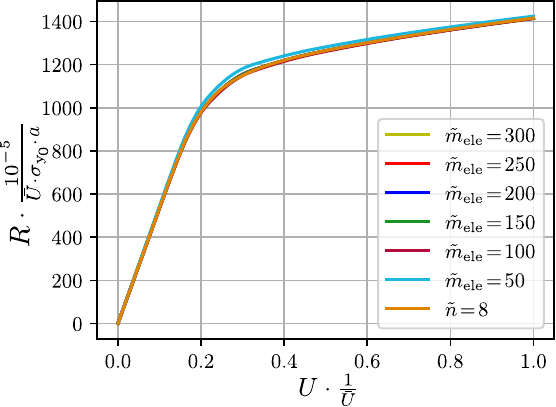}}
	\hfill
	\subfloat[][Relative error w.r.t.~the ROM solution for variation of~$\tilde{m}_\mathrm{ele}$.]{\includegraphics[width=0.47\linewidth]{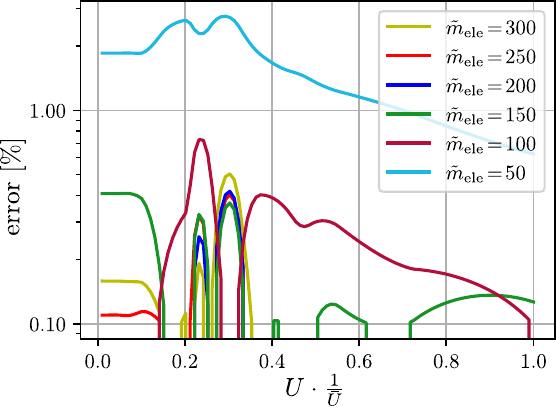}}
	\caption{EHEIM hyper element convergence for a fixed number of $\tilde{n}=8$ ROM modes, when enforcing the criteria of Figure~\ref{fig:DecaySV}\label{fig:EHEIM_conv_study}.}
\end{figure}

\begin{figure}
	\centering
	\includegraphics[width=0.22\linewidth]{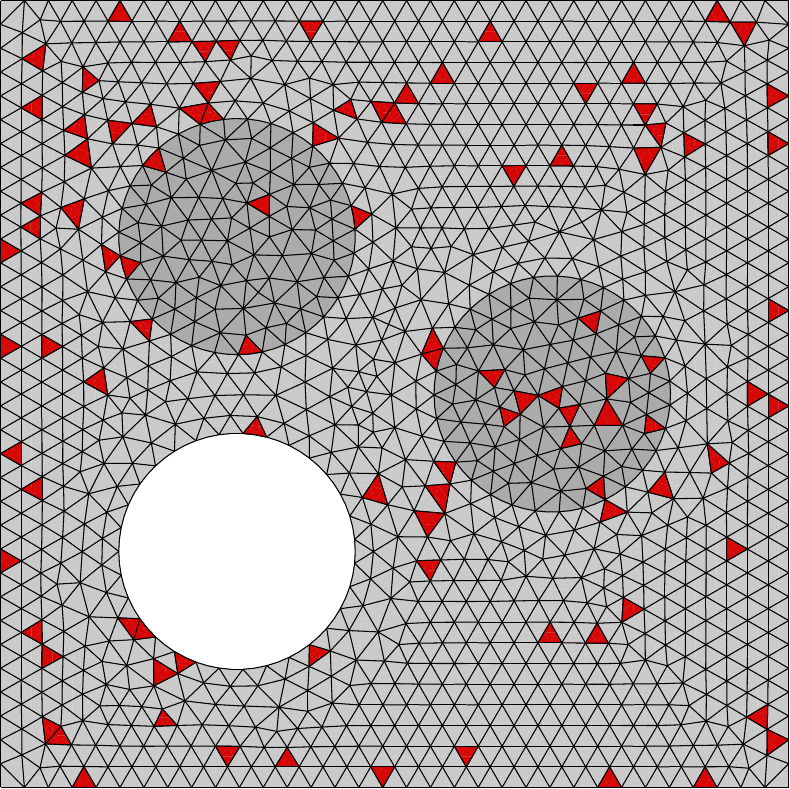}
	\caption{RVE mesh with hyper elements marked in red.}
	\label{fig:Visulization_EHEIM_MieheKoch}
\end{figure}

Finally, the additional criteria are used with the element-based hyperintegration scheme EHEIM, which allows a seamless integration into modular FE frameworks. A convergence study regarding the number of hyper elements is provided in Figure~\ref{fig:EHEIM_conv_study}. 
It is found that $\tilde{m}_\mathrm{ele}=100$ elements are required to reduce the error below 1\%. 
For illustration purposes, the identified hyper elements are highlighted in Figure~\ref{fig:Visulization_EHEIM_MieheKoch}. It is found that they are more or less uniformly distributed over the RVE, with slightly higher concentrations at internal interfaces. Note that only such elements are evaluated throughout the online simulation phase. 
Compared to the corresponding ECM convergence shown in Figure~\ref{fig:EHEIM_conv_study}, the rate of convergence here is slightly inferior. It should be pointed out that each element of the RVE in this example has three integration points. 
However, the overall costs for computing $\underline{\tilde{r}}$ and $\underline{\underline{\tilde{k}}}$ do not scale linearly with the number of integration points in the elements, mainly due to the non-negligible ROM projection costs ($\underline{\underline{\hat{\phi}}}^\mathrm{T}\cdot\bullet$). This observation is further reflected in the respective online costs for the EHEIM simulation in Table~\ref{table:simulation_times_composite} (here, 100 hyper elements, with three integration points each, yield 300 integration points in total). Vice versa, the offline costs of EHEIM are reduced as the training matrix for the necessary costly SDV comprises now only a single row per element, compared to three rows per element with ECM.

\subsection{Hyperelastic porous strip}
\begin{figure}[hbt]
  \begin{minipage}[b]{0.78\textwidth}
    \centering
	\includegraphics[width=0.95\textwidth]{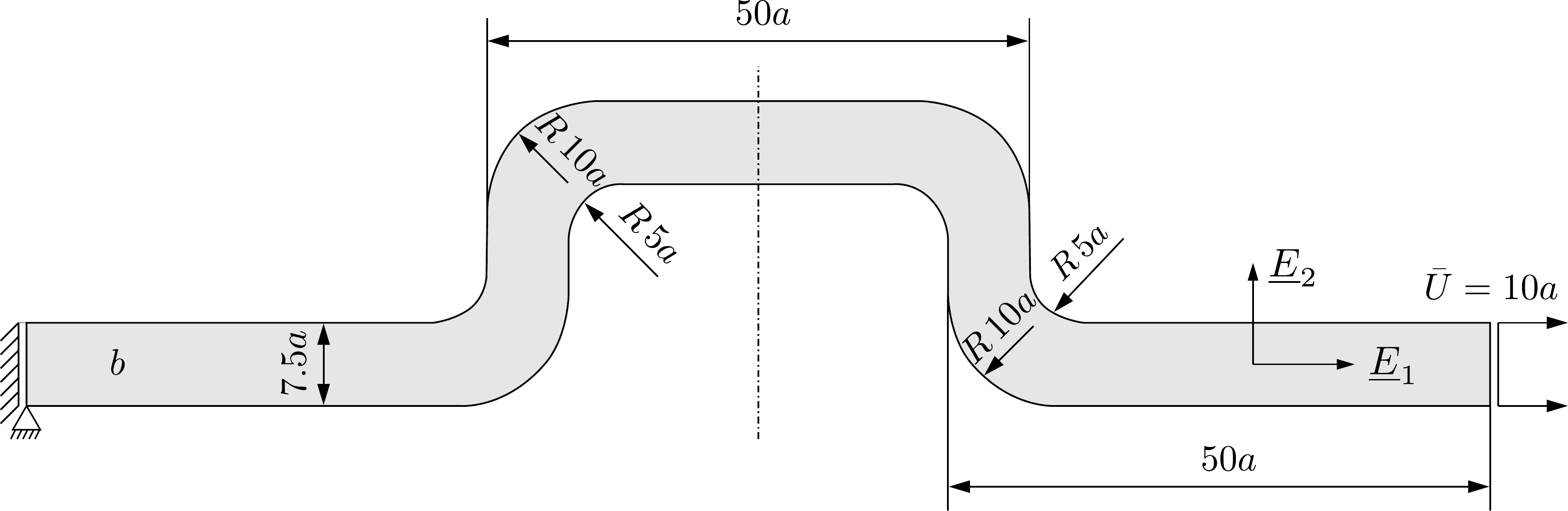}
	\caption{Boundary conditions and dimensions of the macroscopic problem with length parameter $a$ and depth $b$ (where $b=a$).\label{fig:hyperMacro}}      
  \end{minipage}
  \hfill
  \begin{minipage}[b]{0.2\textwidth}
	\includegraphics[width=0.8\textwidth]{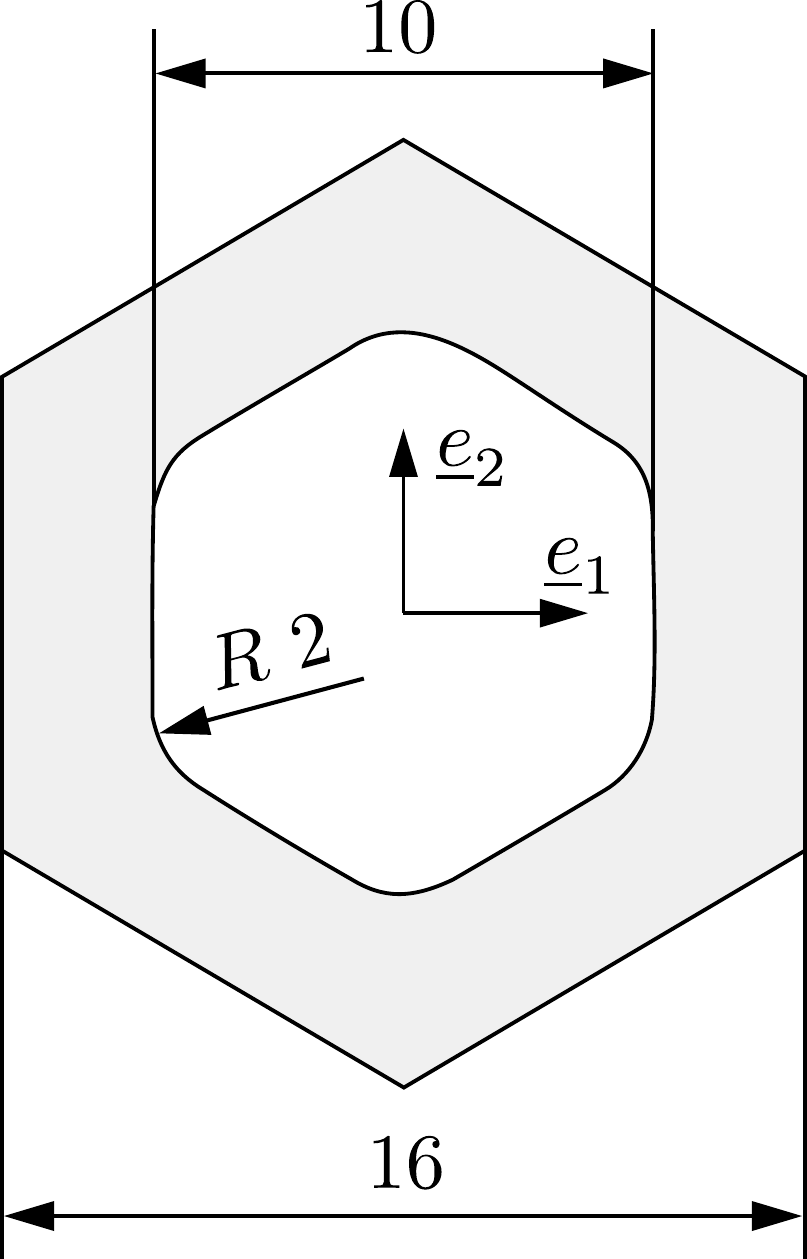}
	\caption{RVE with elastic constants $E,\,\nu=0.3$.\label{RVE_hexa}}
  \end{minipage}
\end{figure}

As a second benchmark problem, a strip that is pulled under displacement control at one side is considered, as shown in Figure~\ref{fig:hyperMacro}, while the observed symmetry will be respected in the simulations. 
The microstructure is defined by a hexagonal RVE with a central pore, see Figure~\ref{RVE_hexa}. Finite deformation theory is applied and the hyperelastic, compressible neo-\textsc{Hooke} law~\cite{Kossa2023}, characterized by the free energy density
\begin{gather}\nonumber
    \psi=\dfrac{\mu}{2}\left[\,\bar{i}_1-3\right]+\dfrac{\kappa}{2}\left[j-1\right]^2\ , \  \\
    \text{with}\ j=\mathrm{det}(\tenq{f})\,,\;\bar{i}_1=\mathrm{tr}(\tenq{\bar{b}})\, ,\; \tenq{\bar{b}}=j^{-2/3}\tenq{f}\cdot\tenq{f}^\mathrm{T}\, ,\;
    \kappa=\dfrac{E}{3\left[1-2\nu\right]}\, , \ \text{and}\ \ \mu=\dfrac{E}{2\left[1+\nu\right]}
\end{gather}
%
is assumed at the micro-scale, for which the implementation described in the \texttt{Abaqus} manual \cite{abaqus} is used in conjunction with a nested plane stress method after \citet{Dodds_1987}.

On both scales, a mixture of rectangular, reduced-integrated and triangular, fully-integrated, quadratic, plane stress standard elements are used. The employed meshes will be shown below,together with the results.

\begin{figure}[hbt]
			\centering
			\includegraphics[width=0.4\textwidth]{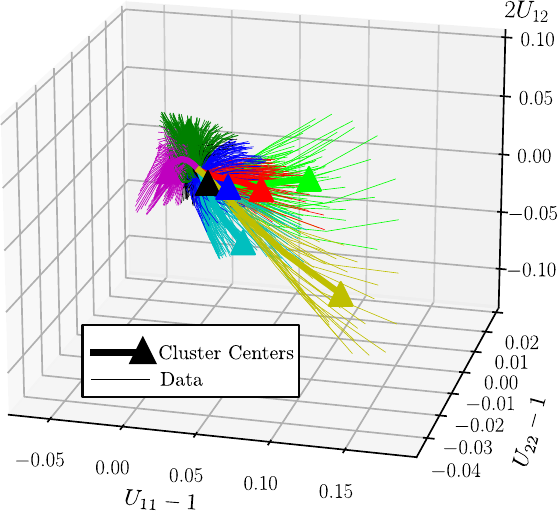}
			\caption{Stretch paths from the simulation with the surrogate model, with eight clusters visualized.}
			\label{fig:Clustering_Visualization}
\end{figure}

\begin{figure}[!hbt]
	{\centering
		\begin{subfigure}[t]{0.48\textwidth}
			\includegraphics[width=0.9\textwidth]{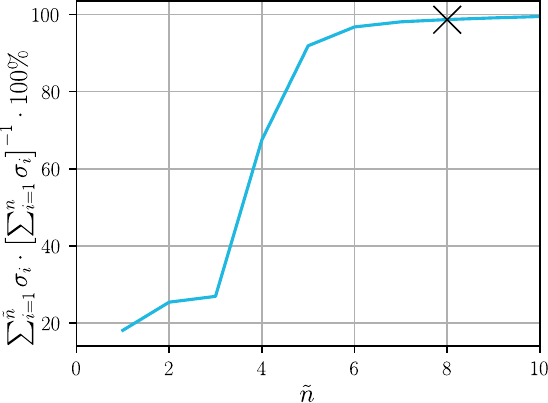}
			\caption{Cumulated singular values of training matrix $\underline{\underline{\hat{x}}}_u$.\label{fig:cumu_sing_xu_wabe}}
		\end{subfigure}\hfill
		\begin{subfigure}[t]{0.48\textwidth}
			\includegraphics[width=0.9\textwidth]{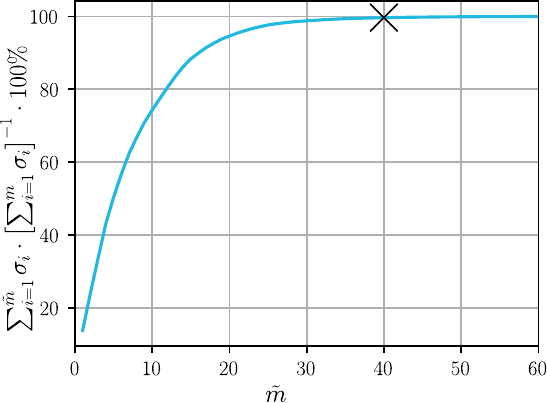}
			\caption{Cumulated singular values of training matrix $\underline{\underline{x}}_f$, where the criteria $\underline{f}_\mathrm{int}$, $\underline{f}_\Sigma$, $\langle\,p\,\rangle$, $\langle\,j\,\rangle$ and $\langle\,\tenq{\sigma}\,\rangle$ were chosen.\label{fig:cumu_sing_xf_wabe}}
		\end{subfigure}\hfill
		\caption{Singular value convergence study, with and without additional selection criteria.}}
\end{figure}

\begin{figure}[!hbt]
			\centering
			\includegraphics[width=0.22\textwidth]{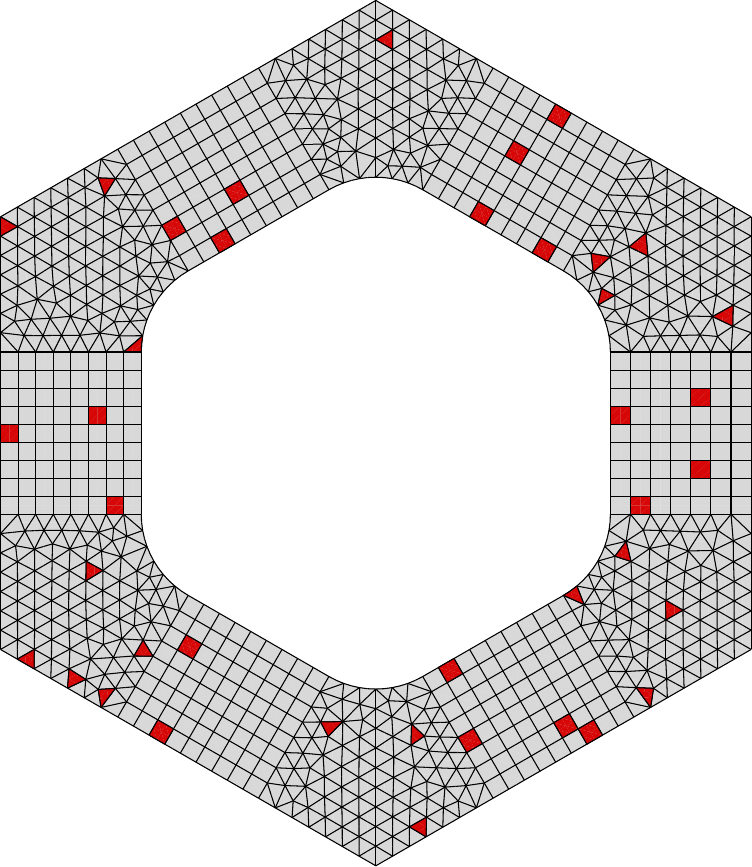}
			\caption{Mesh of the RVE with hyper elements marked in red.}
            \label{fig:Visulization_EHEIM}
\end{figure}

From a simulation of the structure, assuming a homogeneous distribution of the neo-\textsc{Hooke} material, 15 clusters with right stretch trajectories $\tenq{U}(t)$ are obtained. In order to visualize this clustering process, Figure~\ref{fig:Clustering_Visualization} shows the respective paths and eight cluster centers. 
The snapshot matrix $\underline{\underline{\hat{x}}}_u$ is built from HF RVE simulations with these trajectories. In a SVD thereon, a rapid decay of its singular values can be observed, hence their sum converges quickly, as shown in Figure~\ref{fig:cumu_sing_xu_wabe}. The basis matrix $\underline{\underline{\phi}}$ is thus truncated to $\tilde{n}=8$ ROM modes, as indicated in the figure. Subsequently, 40 elements, highlighted in Figure~\ref{fig:Visulization_EHEIM}, are selected by the EHEIM calibration process with the additional criteria, for which the cumulated singular values of $\underline{\underline{x}}_f$ are shown in Figure~\ref{fig:cumu_sing_xf_wabe}, giving evidence that the 40 chosen elements can capture virtually all of the information present in $\underline{\underline{x}}_f$.

\begin{figure}[hbtp]
	\centering
    \includegraphics[width=0.8\textwidth]{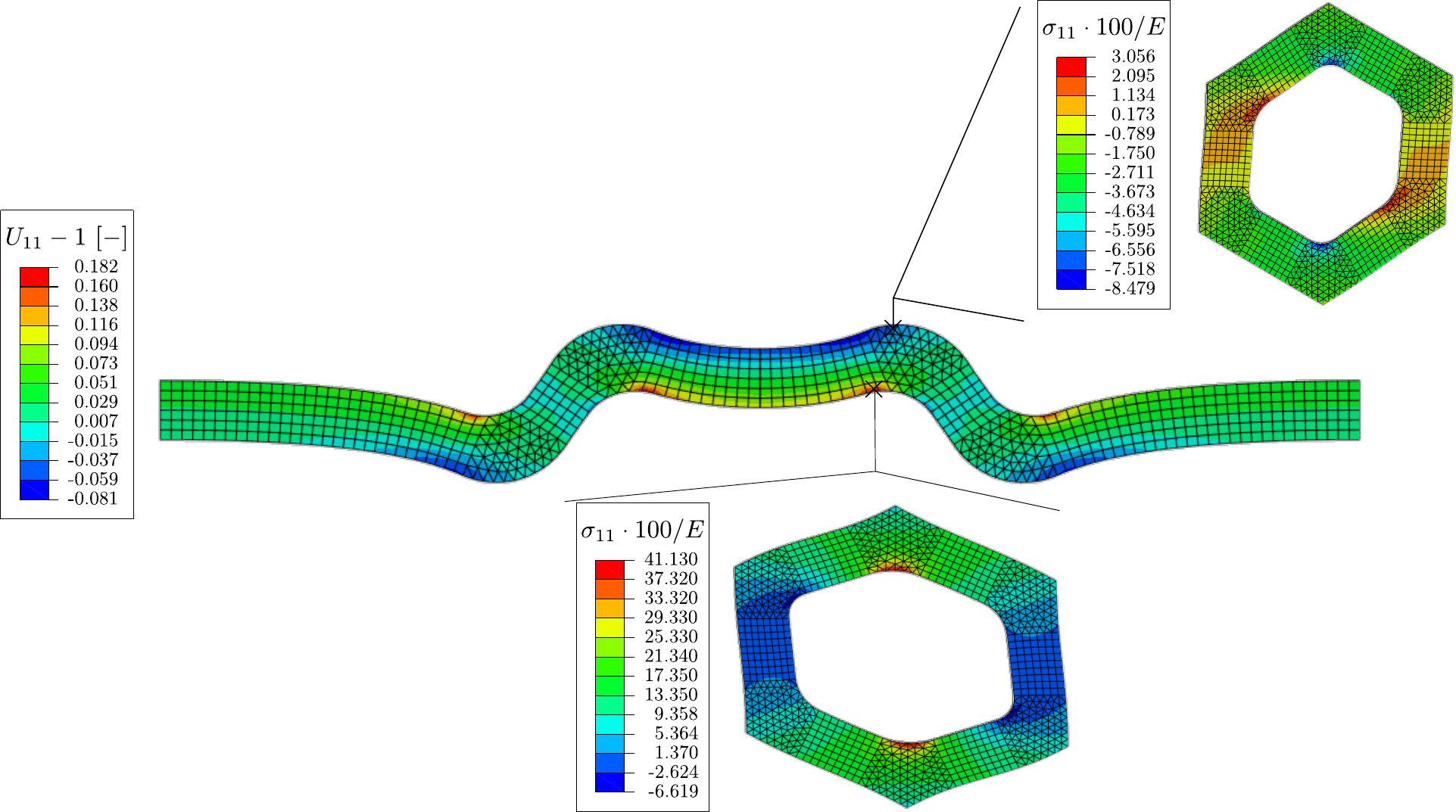}
	\caption{Macroscopic and microscopic fields for the hyperelastic porous strip at maximum load.\label{fig:FEM_results_hyperelastic}}
\end{figure}

\begin{figure}[!hbt]
	\centering
		\begin{subfigure}{0.48\textwidth}
			\includegraphics[width=\textwidth]{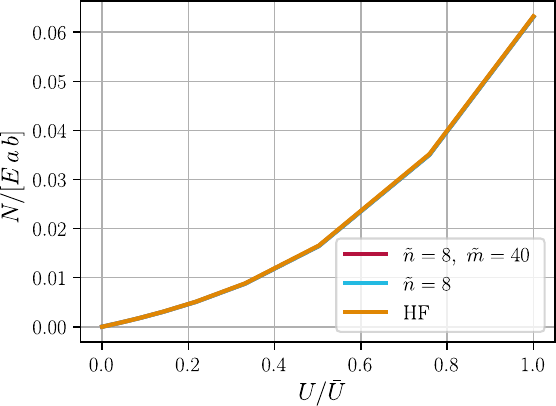}
			\caption{Normalized normal force over prescribed displacement.\label{fig:Norm_force}}
		\end{subfigure}\hfill
		\begin{subfigure}{0.48\textwidth}
			\includegraphics[width=\textwidth]{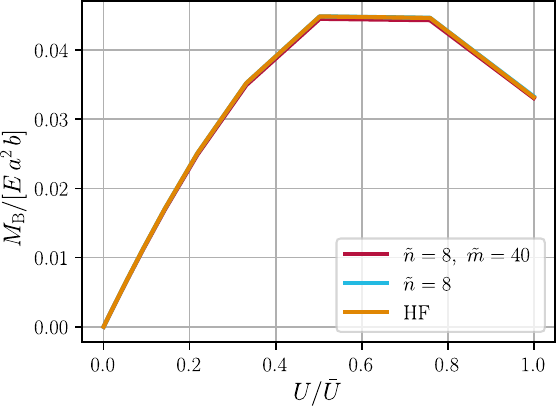}
			\caption{Normalized bending moment over prescribed displacement.\label{fig:BendMom}}
		\end{subfigure}\hfill
		\caption{Global reactions of the hyperelastic porous strip in terms of the applied deflection.}
        \label{fig:results_FM_wabe}
\end{figure}

Detailed FEM simulation results, with inhomogeneous micro-stress and strain distributions, are representatively shown in  Figure~\ref{fig:FEM_results_hyperelastic}. At the strip's end, a normal force ${N=\int_{(A)}\Sigma_{11}\,\mathrm{d}A}$ and bending moment $M_\mathrm{B}=\int_{(A)}\Sigma_{11}\,X_1\,\mathrm{d}A$ appear as a reactions to the prescribed nodal displacements. Their values have been plotted in Figures \ref{fig:Norm_force} and \ref{fig:BendMom} for the different numerical schemes. All methods are seen to agree very well with each other. The nonlinearity in the computed responses is mainly caused by geometrically nonlinear effects. More specifically, when the strip is pulled, it elongates and becomes straighter, which is the reason the normal force $N$ rises progressively and the moment $M_\mathrm{B}$, at first, rises degressively and eventually drops at some point.

\begin{figure}[!hbt]
	\centering
	\includegraphics[width=0.45\textwidth]{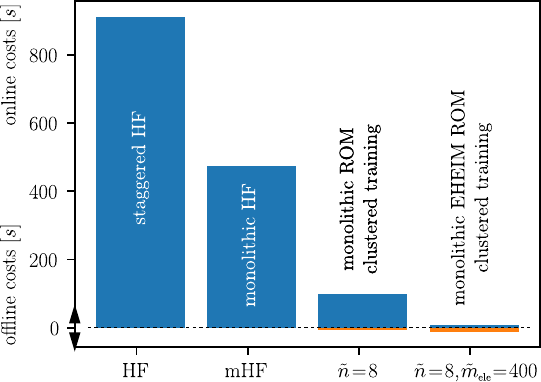}
	\captionof{figure}{Computational online and offline efforts associated with the different numerical schemes for example 2.\label{fig:comp_Effort_wabe}}
\end{figure}

Data on the computational effort for the different methods are compared in Figure~\ref{fig:comp_Effort_wabe} with details given in Table~\ref{table:simulation_times_hyperelastic}. 
\begin{table}[!hbt]
	\centering
	\captionof{table}{Simulation effort for the two-scale, porous, hyperelastic strip problem. The ROM offline time includes simulation of the strip with a surrogate model, RVE-FE training  simulations and SVD on $\underline{\underline{\hat{x}}}_u$; the EHEIM offline time includes ROM-RVE-FE training  simulations, SVD on $\underline{\underline{x}}_f$, and EHEIM selection.\label{table:simulation_times_hyperelastic}}
	\begin{tabular}{|p{0.2\linewidth}|p{0.16\linewidth}|p{0.15\linewidth}|p{0.13\linewidth}|p{0.10\linewidth}|}
		\hline
		method & \multicolumn{2}{c|}{HF} & \multicolumn{2}{c|}{ROM}\\\hline
		SOE & \multicolumn{1}{c|}{staggered} & \multicolumn{3}{c|}{monolithic} \\\hline
		integration & \multicolumn{3}{c|}{full} & \multicolumn{1}{c|}{EHEIM} \\\hline
		DOFs & \multicolumn{2}{c|}{11\,281\,512\,(100\%)} & \multicolumn{2}{c|}{13\,288(0.1\%)}\\\hline
		elements & \multicolumn{3}{c}{2\,393 \,501\,(100\%)} &\multicolumn{1}{|c|}{66\,440\,(2.8\%)}\\\hline
		online time & \multicolumn{1}{c|}{912\,s (100\%)} & \multicolumn{1}{c|}{474\,s (52\%)} & \multicolumn{1}{c}{98\,s (11\%)} & \multicolumn{1}{|c|}{7\,s (0.8\%)}\\\hline
		offline ROM & \multicolumn{2}{c|}{--} & \multicolumn{2}{c|}{7\,s (0.8\%)} \\\hline
		offline EHEIM & \multicolumn{3}{c|}{--}&\multicolumn{1}{c|}{4\,s (0.4\%)}\\\hline
	\end{tabular}
\end{table}
In general, the same trends can be observed as in the first example. More specifically, with an online time of the monolithic EHEIM ROM FE\textsuperscript{2} method below $1\%$ of the conventional FE\textsuperscript{2} method, the speed-up is even slightly higher, while only around another $1\%$ of that time has to be \enquote{invested} offline into the training process. These results demonstrate that the same level of accuracy and speed-up up can be reached at large deformations if the rigid rotation $\tenq{R}$, which could not be represented by an affine projection~\eqref{eq:reducedprojection}, is not transferred to the micro-scale.


\subsection{Woven composite strip with a hole}

\begin{figure}[!hbt]
	\centering
	\begin{subfigure}[b]{0.8\textwidth}
		\includegraphics[width=\textwidth]{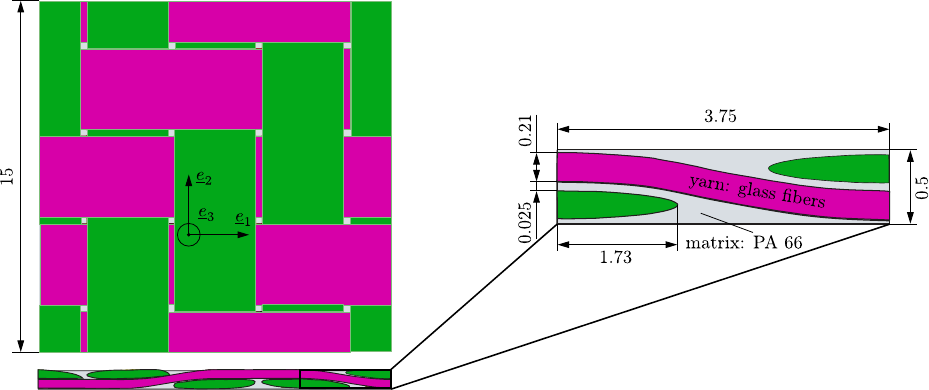}
		\caption{Dimensions of the woven composite meso-structure in [$\mathrm{mm}$].\label{WovenRVE}}
	\end{subfigure}\\[1.5ex]
	\begin{subfigure}[b]{0.55\textwidth}
		\includegraphics[width=\textwidth]{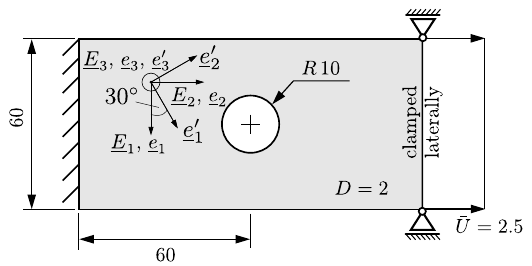}
		\caption{Plate with a  hole problem dimensions (in [$\mathrm{mm}$], depth $D$), fiber orientation and boundary conditions.\label{experimental_setup}}\hfill
	\end{subfigure}\hfill
	\begin{subfigure}[b]{0.4\textwidth}
		\includegraphics[width=\textwidth]{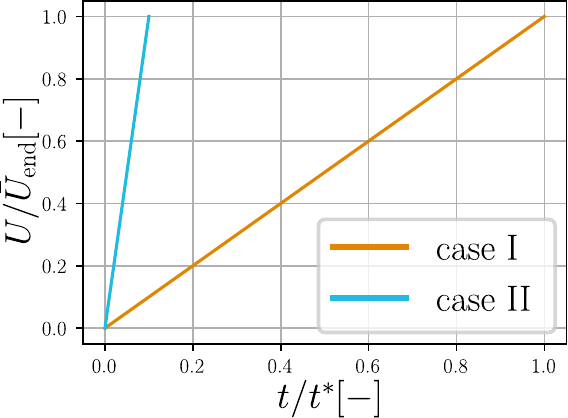}
		\caption{Prescribed displacement over time, for the two considered cases.}
        \label{loadcases}
	\end{subfigure}
	\caption{Microstructure, experimental setup, and loading of a plate with a hole experiment of a woven composite, with two different fiber orientations (after \citet{Tikarrouchine2021}).}
\end{figure}

As a final application, a composite structure of polyamide 6–6 (PA66) reinforced by 2–2 twill weave glass woven fabric is considered, whose behavior was experimentally investigated by \citet{Tikarrouchine2021}.
Its periodic meso-structure is shown in Figure~\ref{WovenRVE}. 
The weft and warp yarns are made of $85\%$ glass fibers embedded in $15\%$ PA 66 matrix. 
The constitutive models of \citet{Praud_2018_dis}, as well as the parameters identified therein, are used for both constituents.
In short, the behavior of the matrix is described by a viscoelastic-viscoplastic law with isotropic damage. 
The yarns, themselves representing a composite of glass fibers within a polymer, are described as a uniaxial composite with micro-cracks, where the effect of the latter is incorporated by means of a Mori-Tanaka micromechanics scheme.

In the original work by Praud et al.~\cite{Praud_2017_matrix,Praud_2018_dis}, a convex cutting plane algorithm was derived for the numerical implementation of this model. However, this approach does not deliver an  algorithmically consistent tangent and can therefore lead to convergence issues, particularly in multi-dimensional loading scenarios. Moreover, in our own experience \cite{Lange_hyper_ROM_2024}, even a closest-point projection general return algorithm with a consistent algorithmic tangent provides unsatisfactory robustness for this complex constitutive model. An ImplEx integration scheme, originally proposed by \citet{Oliver_2008}, is therefore additionally applied here, in order to handle the implicit time integration with a robust regula falsi algorithm. The ImplEx scheme combines the advantage of the unconditional stability of implicit integration with the robustness and global linearity of explicit time integration. The idea is that firstly the state variables $\mathbbl{h}^{t+1}$ are integrated implicitly as usual, and their rates and values are stored for the next time step. Subsequently, updated state variables $\tilde{\mathbbl{h}}^{t+1}$, marked with $(\tilde{\ })$, are computed by an Euler forward step based on the implicitly integrated value $\mathbbl{h}^{t}$ from the last time step. Finally, the values $\tilde{\mathbbl{h}}^{t+1}$ are used within the state law to compute the stress $\tenq{\tilde{\sigma}}^{t+1}$ to be returned to the element routine.
This method is used for both the material models of the polymer matrix and the yarns. The key equations and further information on the constitutive models and their numerical implementation can be found in \ref{sec:micromaterialmodels}.

The macroscopic setup to be simulated is schematically depicted in Figure~\ref{experimental_setup} and comprises a woven composite strip with a circular hole, clamped into the holders of a tensile testing machine. As in the respective experiments \cite{Tikarrouchine2021},  loading of the strip in the FEA is performed in a displacement controlled way, by fixing one end and pulling the other in horizontal direction by $\underline{\bar{U}}(t)$, up to a maximum load $\bar{U}_\mathrm{end}$ at $t^*$.

Four experimental variations are considered in the FEA. Firstly, two RVE orientations are used, as indicated in \figurename~\ref{experimental_setup}, one where the RVE edges are aligned with the axis of loading and one in which the RVE is rotated by $30^\circ$ against the axis of loading. 
Then, two different loading rates $\underline{\dot{\bar{U}}}$ are applied for each orientation, see \figurename~\ref{loadcases}.


In the $0^\circ$ RVE orientation case, a quarter model is sufficient, because material, geometry and loading are symmetric w.r.t.\ to the $x_1$ and $x_2$ coordinate axes at the center of the hole. This material symmetry, of course,  no longer holds in the $30^\circ$ orientation case, so that a full model is required. Macro-scale mesh and time step size convergence studies were performed (not shown here),  using 20-node hexahedral elements with eight integration points (C3D20R in the \texttt{Abaqus} element library). On the micro-scale, the RVE of \citet{Tikarrouchine2021} was used, with standard linear  tetrahedral elements. It is to be emphasized that two-scale FEA for problems of this complexity would hardly have been possible with conventional FE$^2$, due to the computational costs. However, with the EHEIM approach, they can be performed quite easily.  

\begin{figure}
	\begin{subfigure}{\textwidth}
		\includegraphics[width=\textwidth]{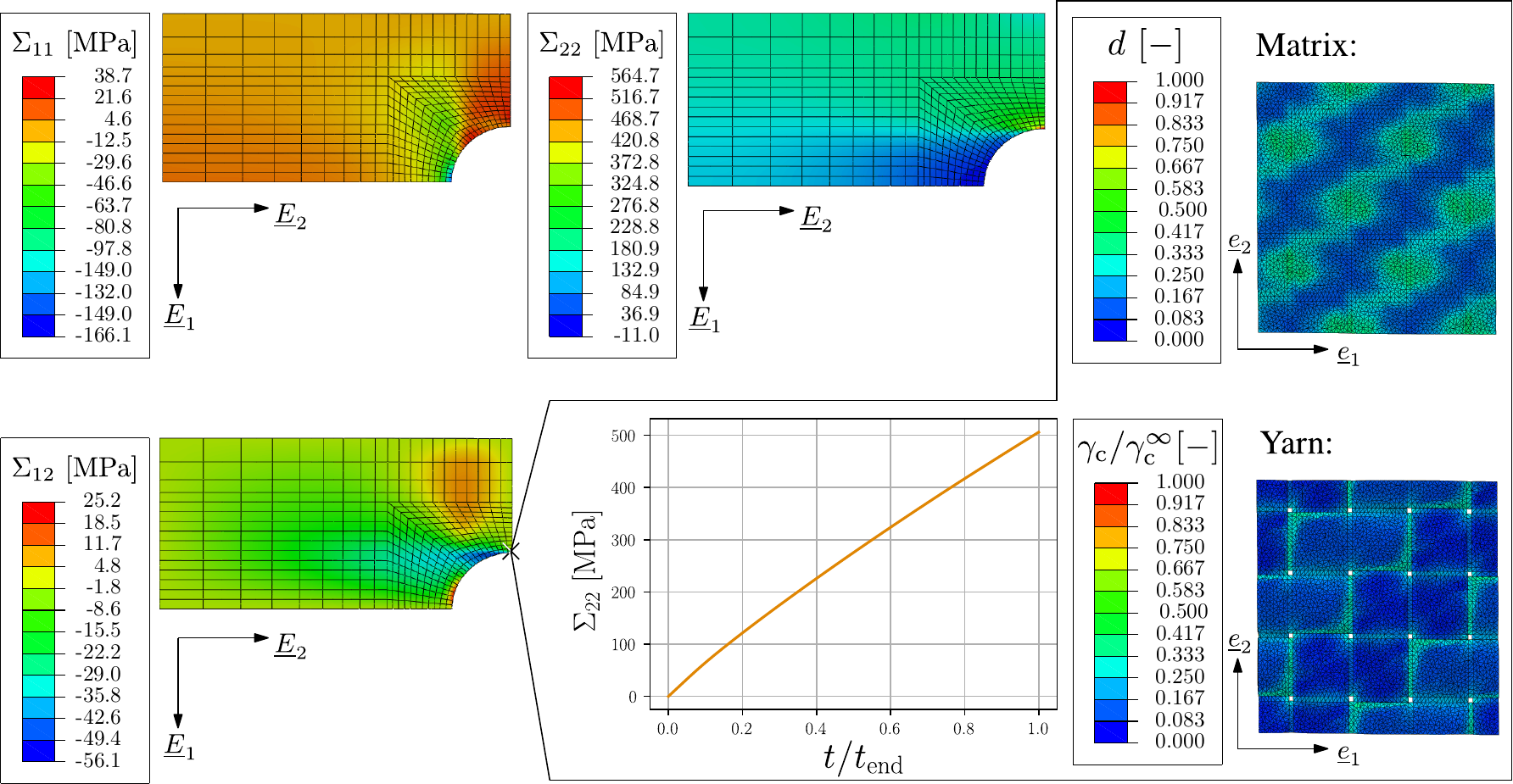}
		\caption{$t=20.83\,\mathrm{s}$, $0^\circ$ RVE orientation, quarter model.\label{results_false_color_0}}
	\end{subfigure}\\[1.0ex]
	\begin{subfigure}{\textwidth}
		\includegraphics[width=\textwidth]{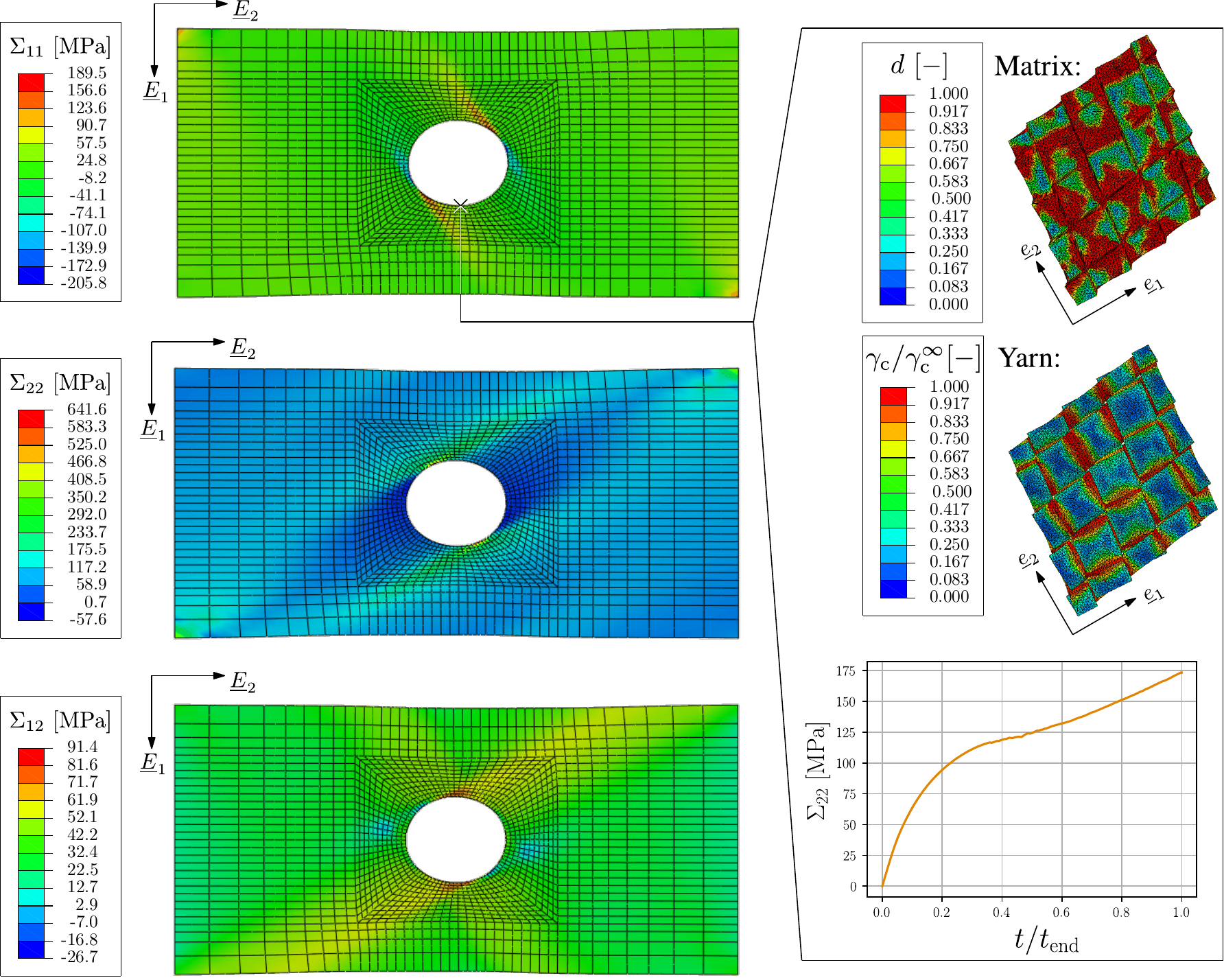}
		\caption{$t=37.88\,\mathrm{s}$, $30^\circ$ RVE orientation, full model.\label{results_false_color_30}}
	\end{subfigure}
	\caption{Detailed numerical results for case I for (a) $0^\circ$ and (b) $30^\circ$ orientation of fibres: (left)  macroscopic stress fields, (right top) damage $d$ in the matrix and micro-crack density $\gamma_{\mathrm\mathrm{c}}$ and (bottom right) macroscopic stress over time $\Sigma_{22}(t)$ for a selected material point \label{detailed:WovenComp}.}
\end{figure}

The analysis, in fact, revealed that a considerably finer mesh, see Figure~\ref{detailed:WovenComp},
is required to resolve the stress concentration near the hole, than was used---or probably even feasible---in the original work \cite{Tikarrouchine2021}. 
Based on the time step convergence study, whose details can be found in the dissertation \cite{LangeDiss}, a step size of $\Delta t=0.01\,t^*$ was chosen for the ImplEx scheme. 
Interestingly, the mesh and time step refinement showed the tendency towards an improved agreement with the experimental results in the global response---although a broader spectrum of numerical studies would have to be performed to truly validate this trend. 

\begin{figure}[!h]
	{\centering
		\begin{subfigure}{0.44\textwidth}
			\includegraphics[width=\textwidth]{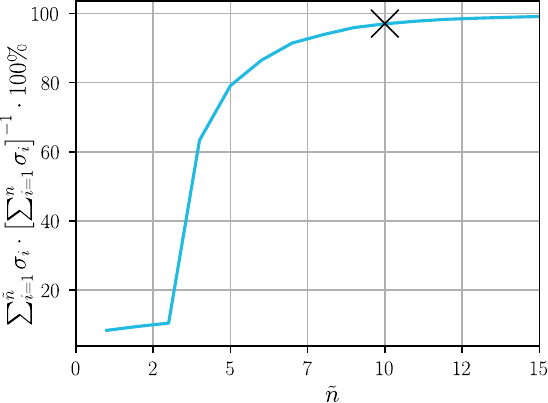}
			\caption{Cumulated singular values of matrix $\underline{\underline{\hat{x}}}_u$.\label{decay_x_u_wovencomp}}
            \vspace{2ex}
		\end{subfigure}\hfill
		\begin{subfigure}{0.44\textwidth}
			\includegraphics[width=\textwidth]{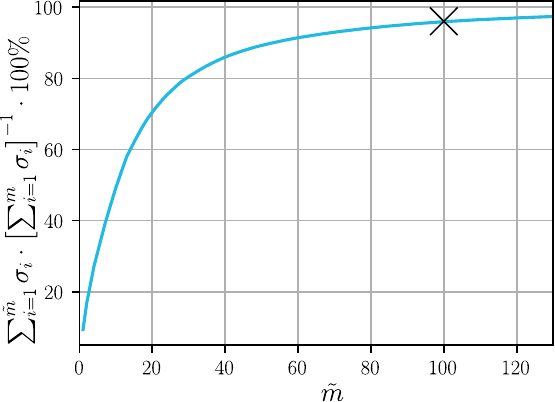}
			\caption{Cumulated singular values of matrix $\underline{\underline{x}}_f$, where the criteria $\underline{f}_\mathrm{int}$, $\underline{f}_\Sigma$, $\langle\,p\,\rangle$ and $\langle\,\tenq{\sigma}\,\rangle$ are chosen.\label{decay_x_f_wovencomp}}
		\end{subfigure}
		\caption{Singular value convergence, , with and without additional selection criteria.}}\vspace{-1ex}
\end{figure}

To employ the clustered hyper ROM method, the macro-structure of the $30^\circ$ orientation was simulated with  linear-elastic properties of the RVE. Here, only the response to the displacement loading rate case II was simulated. By using the clustering method of the strain state, 14 training strain clusters were found. The resulting cumulated singular values of the training matrices $\underline{\underline{\hat{x}}}_u$ and $\underline{\underline{x}}_f$ are plotted in Figures \ref{decay_x_u_wovencomp} and \ref{decay_x_f_wovencomp}, respectively. 
As depicted, ten ROM modes $\underline{\phi}_{\;\!i}$\vspace{0.1ex}, that contain over $95\%$ of the training information and respectively 100 hyper elements, were then chosen in the EHEIM selection process.

\begin{figure}[!h]
	\centering
	\begin{subfigure}{0.44\textwidth}
		\includegraphics[width=\textwidth]{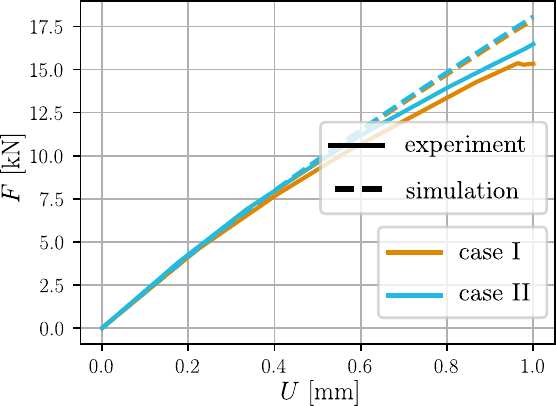}
		\caption{$0^\circ$ RVE orientation, $t^*=20.83\:\mathrm{s}$, $\bar{U}_\mathrm{end}=1.0\:\mathrm{mm}$.\label{fig:curve_wovencomp_0}}
	\end{subfigure}\hfill
	\begin{subfigure}{0.44\textwidth}
		\includegraphics[width=\textwidth]{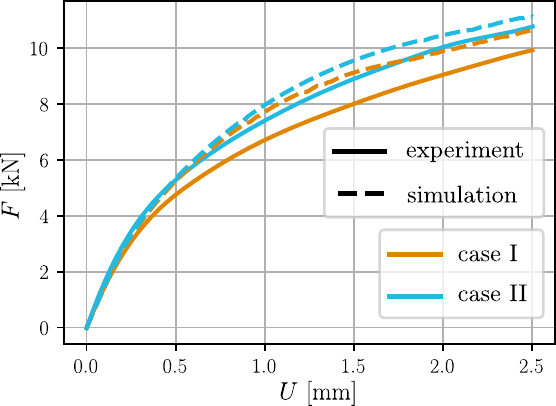}
		\caption{$30^\circ$ RVE orientation, $t^*=37.88\;\mathrm{s}$, $\bar{U}_\mathrm{end}=2.5\:\mathrm{mm}$.\label{fig:curve_wovencomp_30}}
	\end{subfigure}
	\caption{Comparison with the experimental results reported in \citet{Tikarrouchine2021}.\label{fig:curves_wovencomp}}
\end{figure}

The results of the EHEIM-FE$^2$ simulations, in terms of  the computed global force-displacement response curves, are shown in \figurename~\ref{fig:curves_wovencomp}, where they are also compared to the corresponding experimental data \cite{Tikarrouchine2021}. Moreover,  Figure~\ref{detailed:WovenComp} displays detailed local stress and internal state variable distributions on both scales at the point of maximum loading.\footnote{To restore the microscopic fields of the internal state variables, single HF-RVE simulations were performed, with prescribed macroscopic loading history taken from the macroscopic Gauss points of the EHEIM simulations.}
In \figurename~\ref{fig:curve_wovencomp_0}, it can firstly be observed that the component with $0^\circ$ orientation yields a weakly nonlinear behavior, since the load is mainly carried by the yarns. In contrast, the specimens with $30^\circ$ yarn orientation exhibit a more distinct nonlinear load-displacement behavior, see \figurename~\ref{fig:curve_wovencomp_30}. This is explained by the fact that shear stresses  appear, even away from the hole, which have to be transferred by the polymer matrix. Consequently, damage evolves in the matrix, but also in the yarns, due to their relative rotation, as can be seen in \figurename~\ref{results_false_color_30}. At this point, it shall again be emphasized that (only) the employed ImplEx method allows for a robust simulation of such stages of material degradation. 
In both orientations, the faster loading case II leads to higher stresses, and thus global forces, than the slower loading case I, as expected due to the viscous behavior of the material. The relevant experimental results, that are also incorporated in Figure~\ref{fig:curves_wovencomp}, are generally quite well matched by the two-scale FEA, considering the complexity of the problem. Reasons for the remaining discrepancies and potential improvements of the microscopic model are discussed in \cite{Tikarrouchine2021}, but lie beyond the scope of the present work.

\begin{figure}[!h]
	\centering
    \hfill
	\begin{subfigure}[t]{0.38\textwidth}
		\includegraphics[width=\textwidth]{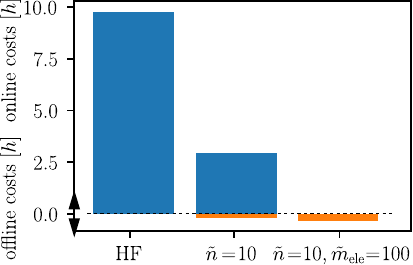}
		\caption{Online and offline simulation times for the three methods: HF, ROM, EHEIM ($0^\circ$ stacking case II, quarter model, coarse mesh).\label{fig:SimTimesActualWoven}}
	\end{subfigure}\hfill
	\begin{subfigure}[t]{0.38\textwidth}
		\includegraphics[width=\textwidth]{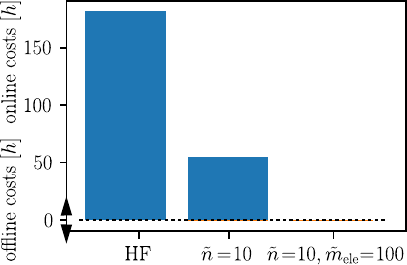}
		\caption{Online and offline simulation times for the EHEIM, extrapolated for HF and ROM  ($30^\circ$ stacking, case II, full model, fine mesh).\label{fig:SimTimesFictiousWoven}}
	\end{subfigure}
    \hfill
	\caption{Computational effort associated with the different numerical schemes.}
\end{figure}

Instead, the present focus is again placed on carefully analyzing the computational costs. To quantify the relative difference of the methods for the woven composite problem, a coarse macro-mesh of the $0^\circ$ orientation simulation (quarter model), with only 64 elements was performed with the HF, ROM and EHEIM methods. This very coarse mesh was chosen so that the HF simulations could finish within reasonable time.
The obtained simulation durations are visualized in Figure~\ref{fig:SimTimesActualWoven}, and listed in Table~\ref{tab:OnlineOfflineTimesWoven}, together with other insightful information. 
\begin{table}[!h]
	\centering
	\captionof{table}{Simulation effort for the woven composite plate with hole problem ($0^\circ$ stacking, case II, quarter model, coarse mesh). The ROM offline time includes simulation of the strip with a surrogate model, RVE-FE training simulations and SVD on $\underline{\underline{\hat{x}}}_u$; the EHEIM offline time includes ROM-RVE-FE training  simulations, SVD on $\underline{\underline{x}}_f$, and EHEIM selection.\label{tab:OnlineOfflineTimesWoven}}
	\begin{tabular}{|p{0.19\linewidth}|p{0.10\linewidth}|p{0.11\linewidth}|p{0.09\linewidth}|}
		\hline
		method & \multicolumn{1}{c|}{HF} & \multicolumn{2}{c|}{ROM}\\\hline
		integration & \multicolumn{2}{c|}{full} & \multicolumn{1}{c|}{EHEIM} \\\hline
		DOFs & \multicolumn{1}{c|}{28\,136\,448\,(100\%)} & \multicolumn{2}{c|}{5\,120\,(0.02\%)}\\\hline
		elements & \multicolumn{2}{c}{55\,905\,280\,(100\%)} &\multicolumn{1}{|c|}{51\,200\,(0.1\%)}\\\hline
		online time & \multicolumn{1}{c|}{35\,213\,s (100\%)} &  \multicolumn{1}{c}{10\,609\,s (30.1\%)} & \multicolumn{1}{|c|}{22\,s (0.06\%)}\\\hline
		offline ROM & \multicolumn{1}{c|}{--} & \multicolumn{2}{c|}{759\,s (2.2\%)} \\\hline
		offline EHEIM & \multicolumn{2}{c|}{--}&\multicolumn{1}{c|}{404\,s (1.1\%)}\\\hline
	\end{tabular}
\end{table}
The HF simulation with the coarse mesh took around 10~h, whereas the EHEIM ROM simulation finished in 22~s, corresponding to a speedup factor of around 1600 regarding the online time and of about 30,  if the total costs (online and offline) are considered.  Moreover, 
the fine mesh, for which symmetry could not be exploited (see Figure~\ref{detailed:WovenComp}, $30^\circ$ stacking), was simulated within 409~s via the EHEIM ROM method---at identical training effort. Hence, when the HF and ROM online simulation times are extrapolated based on the coarse mesh data, as shown in Figure~\ref{fig:SimTimesFictiousWoven}, an overall speedup factor---comprising all offline and online effort---of around 400 was obtained.
Based on these examples, it can be concluded that sufficiently accurate and robust two-scale FEA---even of realistic woven composite structures, accounting for finite deformations and complex, nonlinear, inelastic material behavior---are made possible by the efficiency of the EHEIM ROM method, for which conventional HF FE\textsuperscript{2} simulations would simply not be feasible.

\section{Summary and Outlook}

Hyper reduction, i.e., the efficient computation of the domain integrals appearing in the projected nodal residual, is the key for fast reduced-order concurrent multi-scale simulations. Among the hyper reduction methods discussed in the literature, a quite robust approach turned out to be approximating the projected nodal residual directly from a linear combination of the original Gauss point contributions or element nodal forces. In such an approach, the particular Gauss points or elements to be incorporated, and their respective weights, are determined from the minimization of the error within a set of training simulations. Different measures of the error have been proposed in the literature, such as directly the error within the projected nodal forces in the Empirical Cubature Method (ECM),  or the error within the elastic energy, in the reduced energy optimal cubature (REOC), among others.

The present article combines those methods into a concept of unified integration criteria and proposes the incorporation of additional physically motivated integration criteria. The criteria are implemented for Gauss point-based as well as element-based hyper reduction methods. The proposed method was successfully implemented into the open-source software \texttt{MonolithFE\textsuperscript{2}} \cite{Lange_hyper_ROM_2024,Lange_2021}, which allowed utilizing a previously developed cluster sampling strategy and a monolithic solution of the equilibrium conditions of both scales.

Three different numerical examples were presented, which demonstrated the accuracy, efficiency, and flexibility of the approach. For the complex example of a woven composite with damage at the micro-scale, a robust numerical treatment of the involved material formulations with the aid of the ImplEx time discretization method was additionally introduced.

It turned out that the additional integration criteria improve the accuracy of a hyper-reduced model with given number of hyper integration points, or elements, significantly. Inversely, a considerably lower number of such points, or elements, is necessary to reach a given level of accuracy, so that the incorporation of the proposed unified criteria represents a significant step forward in reducing the computational costs.

Regarding the comparison of Gauss point-based and element-based hyper reduction, it was found that the former allows for a higher online speed-up. The element-based hyper reduction, termed \emph{Empirical Hyper Element Integration Method (EHEIM)}, however, has slightly lower offline costs and, more importantly, allows for a seamless integration into established modular FE frameworks solely at the solver side, making it particularly appealing for flexible multi-scale, multi-physics applications. 

\section*{Acknowledgment}
The authors gratefully acknowledge computing time on the Compute Cluster of the Faculty of Mathematics and Computer Science of Technische Universität Bergakademie Freiberg, operated by the University Computing Center (URZ) and funded by the Deutsche Forschungsgemeinschaft (DFG) under DFG grant number 397252409.

\bibliographystyle{unsrtnat}
\bibliography{bibitems}

\appendix
\section{The Microscale Material Model for Woven Composites}
\label{sec:micromaterialmodels}

\subsection{Viscoelastic-viscoplastic damage of a thermoplastic polymer}\label{chap:matrix_material}
\subsubsection*{Continuous formulation}
Semi-crystalline polymers are used as matrix material in woven composites. In order to describe their complex, rate-dependent behavior, P\textsc{raud} et al.~developed a viscoelastic-viscoplastic-damage constitutive model in the small strain setting, which accurately describes the effects that are observed in experiments with Polyamide 6-6 \cite{Praud_2017_matrix,Praud_2018_dis}. An intensive discussion on the thermodynamical consistent framework can be taken from the aforementioned works. Here only the essential components of the material formulation are discussed and a robust numerical 
treatment is presented.

\begin{figure}[!h]
	\centering
	\includegraphics[width=0.8\textwidth]{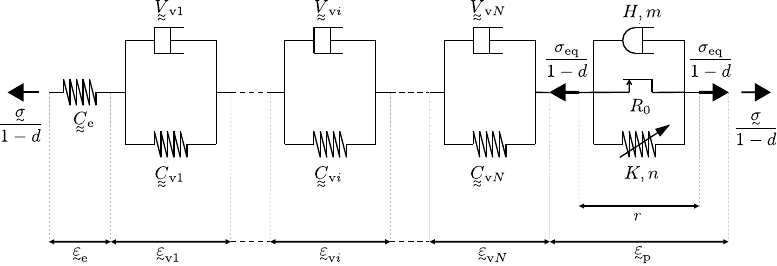}
	\caption{Complete rheological scheme (after P\textsc{raud} \cite{Praud_2018_dis}).}
	\label{fig:RheoScheme}
\end{figure}

The rheological scheme of the material model is shown in Figure~\ref{fig:RheoScheme}. It follows the continuum damage theory with the concept of an effective stress $\tenq{\sigma}/[1-d]$. The model contains an elastic \enquote{spring} rheological element that is characterized by the elasticity tensor $\tenq[2]{C}\phantom{}_\mathrm{e}$ and is associated with the elastic strain $\tenq{\varepsilon}_\mathrm{\!\;e}$. Furthermore, $N$ K\textsc{elvin}-V\textsc{oigt} elements are connected in series to reproduce the viscoelastic behavior at different strain rates, with the viscoelastic strains $\tenq{\varepsilon}_{\mathrm{v}i}$ as internal state variables and corresponding elasticity $\tenq[2]{C}\phantom{}_{\mathrm{v}i}$ and viscosity tensors $\tenq[2]{V}\phantom{}_{\mathrm{v}i}$. The viscoplastic branch incorporates nonlinear hardening, a nonlinear dash-pot, an initial yield stress $R_0$ and is controlled by the internal state variables equivalent viscoplastic strain $r$ and viscoplastic strain~$\tenq{\varepsilon}_\mathrm{\!\;p}$. The evolution of the damage $d$ is coupled to the evolution of the equivalent viscoplastic strain $r$.

The total strain is an additively split of the aforementioned elastic, viscoelastic and viscoplastic strains
\begin{equation}
	\tenq{\varepsilon}=\tenq{\varepsilon}_\mathrm{\!\;e}+\sum_{i=1}^{N}\tenq{\varepsilon}_{\!\;\mathrm{v}i}+\tenq{\varepsilon}_\mathrm{\!\;p}\ 
\end{equation}
Therewith the stress, assuming an isotropic elasticity tensor $\tenq[2]{C}\phantom{}_\mathrm{e}(\kappa,\mu)$ that is represented by aid of the compression $\kappa$ and shear modulus $\mu$ and the spherical $\ProjectionTensSph$ and deviatoric $\ProjectionTensDev$ projection tensors, reads
\begin{equation}\label{state_law_stress}
	\tenq{\sigma}=[\!\;1-d\!\;]\,\tenq[2]{C}\phantom{}_\mathrm{e}:\tenq{\varepsilon}_\mathrm{\!\;e}=[\!\;1-d\!\;]\,\tenq[2]{C}\phantom{}_\mathrm{e}:\left[\tenq{\varepsilon}-\tenq{\varepsilon}_\mathrm{\!\;p}-\sum_{i=1}^{N}\tenq{\varepsilon}_{\!\;\mathrm{v}i}\right]\ ,\ \ \ \tenq[2]{C}\phantom{}_\mathrm{e}=3\,\mathrm{\kappa}\,\ProjectionTensSph+2\,\mathrm{\mu}\,\ProjectionTensDev\ ,\ \kappa=\dfrac{E}{3\left[1-2\nu\right]}\ \ \text{and}\ \ \mu=\dfrac{E}{2\left[1+\nu\right]}.
\end{equation}
The viscoelastic stresses are derived from the elasticity law with the same P\textsc{oisson}'s ratio $\nu$ from above as
\begin{equation}\label{el_law_vi}
	\tenq{\sigma}_{\mathrm{v}i}=\tenq[2]{C}\phantom{}_{\mathrm{v}i}:\tenq{\varepsilon}_{\!\;\mathrm{v}i}\ ,\ \ \ \tenq[2]{C}\phantom{}_{\mathrm{v}i}=3\,\mathrm{\kappa}_{\mathrm{v}i}\,\ProjectionTensSph+2\,\mathrm{\mu}_{\mathrm{v}i}\,\ProjectionTensDev\ .
\end{equation}
The evolution of the viscoelastic strains is given with the isotropic viscosity tensors $\tenq[2]{{V}}\phantom{}_{\mathrm{v}i}$ by
\begin{equation}\label{viso_el_evol}
	\tenq{\dot{\varepsilon}}_{\!\;\mathrm{v}i}\!=\!\tenq[2]{{V}}\phantom{}_{\mathrm{v}i}^{-1}\!:\!\left[\dfrac{\tenq{\sigma}}{1-d}-\tenq{\sigma}_{\mathrm{v}i}\right]\!=\!\dfrac{E}{\eta_{\mathrm{v}i}}\!\left[\tenq{\varepsilon}-\sum_{i=1}^{N}\tenq{\varepsilon}_{\!\;\mathrm{v}i}-\tenq{\varepsilon}_\mathrm{\!\;p}\right]-\dfrac{E_{\mathrm{v}i}}{\eta_{\mathrm{v}i}}\tenq{\varepsilon}_{\!\;\mathrm{v}i}\ ,
\end{equation}
where the fact was used that
\begin{equation}
	\tenq[2]{{V}}\phantom{}_{\mathrm{v}i}^{-1}\!:\tenq[2]{{C}}\phantom{}_\mathrm{e}\!=\!\dfrac{E}{\mathrm{\eta}_{\mathrm{v}i}}\,\tenq[2]{{I}}\ \ \ \text{and}\ \ \ \tenq[2]{{V}}\phantom{}_{\mathrm{v}i}^{-1}\!:\tenq[2]{{C}}\phantom{}_{\mathrm{v}j}\!=\!\dfrac{E_{\mathrm{v}i}}{\mathrm{\eta}_{\mathrm{v}j}}\,\tenq[2]{{I}}\ .
\end{equation}
The equivalent plastic strain $r$ evolves with a certain rate when the yield surface $f$ is exceeded as indicated by the M\textsc{acaulay}-brackets. The evolution is given in terms of the power law
\begin{equation}\label{pl_mult_evol}
	H\,\dot{r}^m=\langle\,f\,\rangle
\end{equation}
with the positive parameters $H$ and $m$. This relation corresponds to the irreversibility of the process i.e.~the equivalent plastic strain $r$ can only grow. The $J_2$-type yield surface is governed by an initial yield stress $R_0>0$ and a power law hardening function $R(r)$ with parameters $K>0$ and $n>0$ and reads
\begin{equation}\label{yield_surf}
	f=\dfrac{\sigma_\mathrm{eq}}{1-d}-R(r)-{R_0} \ \ \ \mathrm{with}\ \ R(r)=K\,r^n\ .
\end{equation}
The evolution of the viscoplastic strain is given by the following associated flow rule
\begin{equation}\label{viso_pl_evol}
	\tenq{\dot{\varepsilon}}_\mathrm{\!\;p}=\dfrac{\partial\,f}{\partial\,\tenq{\sigma}}\dot{r}=\dfrac{3\,\dot{r}}{2\,\sigma_\mathrm{eq}\,\left[1-d\right]}\,\tenq{\sigma}_\mathrm{d}\ \ \ .
\end{equation}
The rate of the damage $\dot{d}$ is coupled with the rate of the equivalent plastic strain $\dot{r}$ over the following evolution equation
\begin{equation}\label{damage_evol}
	\dot{d}=\left[\dfrac{Y}{S}\right]^\mathrm{\beta}\,\dfrac{\dot{r}}{1-d}\ ,\ \ \ 0\le d\le 1.
\end{equation}
In there, $S>0$ and $\beta$ are parameters and $Y$ corresponds to the total energy release density given by
\begin{equation}\label{energie_release}
	Y=\dfrac{1}{2}\ \tenq{\varepsilon}_\mathrm{\!\;e}:\tenq[2]{{C}}\phantom{}_\mathrm{e}:\tenq{\varepsilon}_\mathrm{\!\;e}+\sum_{i=1}^{N}\left(\,\dfrac{1}{2}\ \tenq{\varepsilon}_{\!\;\mathrm{v}i}:\tenq[2]{{C}}\phantom{}_{\mathrm{v}i}:\tenq{\varepsilon}_{\!\;\mathrm{v}i}\right)
\end{equation}
Since it is ensured that $\dot{r}>0$ it follows that the damage can only grow $\dot{d}>0$.
\subsubsection*{Implicit integration}

At first the evolution equations \eqref{viso_el_evol}, \eqref{pl_mult_evol}, \eqref{viso_pl_evol} and \eqref{damage_evol} are discretized with the E\textsc{uler} \emph{backward} method
\begin{equation}\label{discr_vi}
	\tenq{\varepsilon}_{\;\!\mathrm{v}i}^{t+1}=\tenq{\varepsilon}_{\;\!\mathrm{v}i}^{t}+\Delta t\,\dfrac{E}{\eta_{\mathrm{v}i}}\left[\tenq{\varepsilon}^{t+1}-\sum_{i=1}^{N}\tenq{\varepsilon}_{\;\!\mathrm{v}i}^{t+1}-\tenq{\varepsilon}_\mathrm{\;\!p}^{t+1}\right]-\Delta t\,\dfrac{E_{\mathrm{v}i}}{\eta_{\mathrm{v}i}}\,\tenq{\varepsilon}_{\;\!\mathrm{v}i}^{t+1}\ ,
\end{equation}
\begin{equation}\label{discr_r}
	r^{t+1}=r^{t}+\dfrac{\Delta t}{\mathrm{H}^{1/\mathrm{m}}}\,\langle\,f^{t+1}\,\rangle^{1/\mathrm{m}}=r^{t}+\dfrac{\Delta t}{\mathrm{H}^{1/\mathrm{m}}}\left\langle\dfrac{\sigma_\mathrm{eq}^{t+1}}{1-d^{\;\!t+1}}-K\,\left[r^{t+1}\right]^n-{R_0}\right\rangle^{\frac{1}{m}}\ ,
\end{equation}
\begin{equation}\label{discr_ep}
	\tenq{\varepsilon}_\mathrm{\;\!p}^{t+1}=\tenq{\varepsilon}_\mathrm{\;\!p}^{t}+\dfrac{3\,\left[r^{t+1}-r^{t}\right]}{2\,\sigma_\mathrm{eq}^{t+1}\ \left[1-d^{\;\!t+1}\right]}\,\tenq{\sigma}_\mathrm{d}^{t+1}
\end{equation}
and
\begin{equation}\label{discr_d}
	d^{\;\!t+1}=d^{\;\!t}+\left[\dfrac{Y^{t+1}}{S}\right]^\mathrm{\beta}\,\dfrac{r^{t+1}-r^{t}}{1-d^{\;\!t+1}}\ .
\end{equation}
The evolution equation \eqref{discr_d} for the damage is a problem from a numerical standpoint, because the damage cannot be larger than one, and a fully damaged integration point has zero stiffness $\tenq[2]{{C}}_\mathrm{t}=\tenq[2]{{0}}$ and will result in problems when solving the global system of equations. To circumvent this issue, equation \eqref{discr_d} is modified after being put into residual form 
\begin{equation}\label{discr_d_mod}
	R_d\!:=d^{\;\!t+1}\!-d^{\;\!t}\!-\!\left[\dfrac{Y^{t+1}}{S}\right]^\mathrm{\beta}\!\![r^{t+1}-r^{t}]
	\begin{cases}
		\!\dfrac{1}{1-d^{\;\!t+1}}&, d^{\;\!t+1}<0.99\\[2ex]
		\!\left[3.0\left(d^{\;\!t+1}\right)^3\!-8.96\left(d^{\;\!t+1}\right)^2\!+8.92\,d^{\;\!t+1}-2.96\right]\!10^8&, \text{o/w}
	\end{cases}\ ,
\end{equation}
such that the function is sufficiently smooth, and the damage will always grow, but zero damage cannot be reached.

The discretized evolution equations can be inserted into one another, such that a single nonlinear equation for $r^{t+1}$ emerges. The first evolution equations \eqref{discr_vi} form a system of equations for the viscoelastic strains with the viscoplastic strains as unknowns on the right-hand side
\begin{equation}\label{system_of_eq_vi}
	\begin{bmatrix}
		\tenq{\varepsilon}_{\;\!\mathrm{v}1}^{t+1}\\[2ex]
		\tenq{\varepsilon}_{\;\!\mathrm{v}2}^{t+1}\\[2ex]
		...\\[2ex]
		\tenq{\varepsilon}_{\;\!\mathrm{v}\mathrm{N}}^{t+1}
	\end{bmatrix}\!\!=
\underline{\underline{A}}^{-1}\,\cdot
	\begin{bmatrix}
		\tenq{\varepsilon}_{\;\!\mathrm{v}1}^{t}+\frac{\Delta t\,E}{\eta_{\mathrm{v}1}}\,\left[\tenq{\varepsilon}^{t+1}-\tenq{\varepsilon}_\mathrm{\;\!p}^{t+1}\right]\\[2ex]
		\tenq{\varepsilon}_{\;\!\mathrm{v}2}^{t}+\frac{\Delta t\,E}{\eta_{\mathrm{v}2}}\,\left[\tenq{\varepsilon}^{t+1}-\tenq{\varepsilon}_\mathrm{\;\!p}^{t+1}\right]\\[2ex]
		...\\[2ex]
		\tenq{\varepsilon}_{\;\!\mathrm{v}\mathrm{N}}^{t}+\frac{\Delta t\,E}{\eta_{\mathrm{v}\mathrm{N}}}\,\left[\tenq{\varepsilon}^{t+1}-\tenq{\varepsilon}_\mathrm{\;\!p}^{t+1}\right]
	\end{bmatrix}\!,\ 
	A_{ij}=
     \begin{cases}
	1\!+\!\frac{\Delta t\,E}{\eta_{\mathrm{v}i}}\!+\!\frac{\Delta t\,E_{\mathrm{v}i}}{\eta_{\mathrm{v}i}}, & \!\!i=j \\[2ex]
	\frac{\Delta t\,E}{\eta_{\mathrm{v}i}}, & \!\!i\ne j
\end{cases}\ ,\ \ \underline{\underline{A}}\in\mathbb{R}^{N\times N}\!.
\end{equation}
In a trial step it has to be assumed that the yield surface is not exceeded such that the plastic strain does not evolve $r^{t+1}=r^{t}$ from which follows $	\tenq{\varepsilon}_\mathrm{\;\!p}^{t+1}=\tenq{\varepsilon}_\mathrm{\;\!p}^{t}$. Then, the viscoelastic strains $\tenq{\varepsilon}_{\;\!\mathrm{v}i}^{t+1}$ can be computed. If the therewith evaluated yield function \eqref{yield_surf} is smaller than zero the step will be viscoelastic and the trial viscoelastic strains were the actual ones. However, if the yield surface is exceeded, a correction must be performed. In doing so, at first the sum of the viscoelastic strains is rephrased using equation \eqref{system_of_eq_vi}
\begin{equation}\label{res_vi}
	\tenq{\varepsilon}_\mathrm{\;\!v}^{t+1}:=\sum_{i=1}^{N}	\tenq{\varepsilon}_{\mathrm{\;\!v}i}^{t+1}=\underbrace{\sum_{i=1}^{\mathrm{N}}\sum_{j=1}^{\mathrm{N}}A_{ij}^{-1}\,\tenq{\varepsilon}_{\mathrm{\;\!v}j}^{t}}_{\textstyle:=\tenq{\varepsilon}_\mathrm{\;\!v}^{t\,\mathrm{mod}}}+\underbrace{\sum_{i=1}^{\mathrm{N}}\sum_{j=1}^{\mathrm{N}}\left(A_{ij}^{-1}\frac{\Delta t\,E}{\eta_{\mathrm{v}j}}\right)}_{\textstyle:=B}\,\left[\tenq{\varepsilon}^{t+1}-\tenq{\varepsilon}_\mathrm{\;\!p}^{t+1}\right]\,,
\end{equation}
defining the modified viscous strain $\tenq{\varepsilon}_\mathrm{\;\!v}^{t\,\mathrm{mod}}$ and the (precomputable) scalar $B$ as denoted.
This result can be inserted in the elasticity law \eqref{state_law_stress} in order to get the stress deviator
\begin{align}
	\nonumber\tenq{\sigma}_\mathrm{d}^{\;\!t+1}
	&=2\mu\,\left[1-d^{\;\!t+1}\right]\,\left[(1-B)\,\tenq{\varepsilon}_\mathrm{\;\!d}^{t+1}-\tenq{\varepsilon}_\mathrm{\;\!v\,d}^{t\,\mathrm{mod}}-(1-B)\,\tenq{\varepsilon}_\mathrm{\;\!p}^{t+1}\right]\ .
\end{align}
The viscoplastic strain $\tenq{\varepsilon}_\mathrm{\;\!p}^{t+1}$ can be removed by inserting the discretized evolution equation~\eqref{discr_ep}, which gives an expression for the deviatoric stress
\begin{align}\label{stress_dev}
	\tenq{\sigma}_\mathrm{d}^{\;\!t+1}&\!=\!2\mu\left[1\!-d^{\;\!t+1}\right]\!\dfrac{\sigma_\mathrm{eq}^{t+1}}{\sigma_\mathrm{eq}^{t+1}+3\,\mu[1-B]\left[r^{t+1}-r^{t}\right]}\!\left[(1-B)\,\tenq{\varepsilon}_\mathrm{\;\!d}^{t+1}\!-\tenq{\varepsilon}_\mathrm{\;\!v\,d}^{t\,\mathrm{mod}}\!-\!(1-B)\,\tenq{\varepsilon}_\mathrm{\;\!p}^{t}\right].
\end{align}
When the v.~M\textsc{ises} norm is taken from both sides of equation \eqref{stress_dev}
\begin{equation}\label{eq_stress_final}
	\sigma^{t+1}_\mathrm{eq}=2\mu\,\left[1-d^{\;\!t+1}\right]\,\left[(1-B)\,\tenq{\varepsilon}_\mathrm{\;\!d}^{t+1}-\tenq{\varepsilon}_\mathrm{\;\!v\,d}^{t\,\mathrm{mod}}-(1-B)\,\tenq{\varepsilon}_\mathrm{\;\!p}^{t}\right]_\mathrm{eq}-3\mu\,[1-B]\left[r^{t+1}-r^{t}\right]\ ,
\end{equation}
an expression for the equivalent stress $\sigma^{t+1}_\mathrm{eq}$ is gained, which only depends on $d^{\;\!t+1}$ and $r^{t+1}$. Equation \eqref{eq_stress_final} can be inserted into the discretized evolution equation \eqref{discr_r} for $r^{t+1}$ and put into residual form, yielding
\begin{align}\label{eqn:res_r}
	\nonumber R_r\!:=r^{t+1}-\,r^{t}-\dfrac{\Delta t}{\mathrm{H}^{1/\mathrm{m}}}\Bigg[&2\mu\left[(1-B)\,\tenq{\varepsilon}_\mathrm{\;\!d}^{t+1}-\tenq{\varepsilon}_\mathrm{\;\!v\,d}^{t\,\mathrm{mod}}-(1-B)\,\tenq{\varepsilon}_\mathrm{\;\!p}^{t}\right]_\mathrm{eq}\\
	&-\dfrac{3\mu[1-B]\!\left[r^{t+1}\!\!-r^{t}\right]}{1-d^{\;\!t+1}}-K\!\left[r^{t+1}\right]^n-\mathrm{R_0}\Bigg]^{\!\frac{1}{m}}\ .
\end{align}
Now the regula falsi method, which does always converge, is applied to solve $R_d(d^{\;\!t+1})=0$ for $d^{\;\!t+1}$ and in a staggered manner with the Newton-Raphson method to solve $R_r(d^{\;\!t+1},r^{t+1})=0$ for $r^{t+1}$. 
The value from the last time increment $d^{\;\!a}=d^{\;\!t}$ is taken is first initial guess. The required second initial guess value $d^{\;\!b}$ is found by interval interval nesting until $R_{\;\!d}^{\;\!b}$ has a different sign than $R_{\;\!d}^{\;\!a}$. 

If values $d^{\;\!t+1}$ and $r^{t+1}$ are found, $\sigma_\mathrm{eq}^{t+1}$, $\tenq{\sigma}^{t+1}_\mathrm{d}$, $\tenq{\varepsilon}_\mathrm{\;\!p}^{t+1}$, $\tenq{\varepsilon}_{\;\!\mathrm{v}i}^{t+1}$, $\tenq{\sigma}_{\mathrm{v}i}^{t+1}$, $\tenq{\sigma}^{t+1}$ and $Y^{t+1}$ can be updated via equations \eqref{eq_stress_final}, \eqref{stress_dev}, \eqref{discr_ep}, \eqref{system_of_eq_vi}, \eqref{el_law_vi}, \eqref{state_law_stress} and  \eqref{energie_release}, respectively.

\subsubsection*{Explicit integration for ImplEx}
The explicit measures of the material formulation, under which are the stress and algorithmic tangent that will be returned to the element routine, are marked with a tilde $(\tilde{\phantom{a}})$. The measures from the last step without tilde are those computed with the E\textsc{uler} \emph{backward} method. In order to evaluate the stress state law \eqref{state_law_stress}, it is necessary to know $\tilde{d}^{\,t+1}$ and  $\tenq{\tilde{\varepsilon}}_\mathrm{\;\!p}^{\,t+1}$. The damage is simply explicitly integrated using the rate of the last step
\begin{equation}
	\tilde{d}^{\,t+1}=d^{\,t}+\Delta t\,\dot{d}^{\,t}=d^{\,t}+\Delta t\,\dfrac{d^{\,t}-d^{\,t-1}}{t^{t}-t^{t-1}}\ .
\end{equation}
This kind of explicit integration is not a suitable option for $\tenq{\tilde{\varepsilon}}_\mathrm{\;\!p}^{\,t+1}$. In the work of O\textsc{liver} et al.~\cite{Oliver_2008} it was shown, that rephrasing the plastic flow as shown below yields a very cheap way to compute the stress and obtain a positive definite, isotropic algorithmic tangent. The rephrased plastic flow reads
\begin{equation}\label{exp_damage}
	\dot{\bar{\lambda}}=\dfrac{3}{2\,\sigma_\mathrm{eq}\,[1-d]}\,\dot{r}\ \ \ \text{and}\ \ \ \tenq{\dot{\varepsilon}}_\mathrm{\;\!p}=\dot{\bar{\lambda}}\ \tenq{\sigma}_\mathrm{d}
\end{equation}
which is just an alternative formulation of equation \eqref{viso_pl_evol} with a differently defined plastic multiplier $\dot{\bar{\lambda}}$. In discretized form, the explicit integration becomes therewith
\begin{equation}\label{exp_vipl}
	\Delta\tilde{\bar{\lambda}}^{t+1}=\dfrac{3}{2\,\sigma_\mathrm{eq}^{t}\,[1-d^{\;\!t}]}\,\dfrac{r^t-r^{t-1}}{t^{t}-t^{t-1}}\,\Delta t\ \ \ \text{and}\ \ \ \tenq{\tilde{\varepsilon}}_\mathrm{\;\!p}^{\;\!t+1}=\tenq{\varepsilon}_\mathrm{\;\!p}^{t}+\Delta\tilde{\bar{\lambda}}^{t+1}\,\tenq{\tilde{\sigma}}_\mathrm{d}^{t+1}=\tenq{\varepsilon}_\mathrm{\;\!p}^{t}+\Delta\tilde{\bar{\lambda}}^{t+1}\,\ProjectionTensDev:\tenq{\tilde{\sigma}}^{t+1}\ .
\end{equation}
With aid of equation \eqref{res_vi} and the obtained results \eqref{exp_damage} and \eqref{exp_vipl}, the elasticity law~\eqref{state_law_stress} can be evaluated, which yields
\begin{equation}
	\tenq{\tilde{\sigma}}^{t+1}\!=\!(1-\tilde{d}^{\;\!t+1})\!\left[3\,\mathrm{\kappa}\,\ProjectionTensSph+2\,\mathrm{\mu}\,\ProjectionTensDev\right]\!\!:\!\!\left[(1\!-\!B)\,\tenq{\varepsilon}^{t+1}\!-\tenq{\varepsilon}_\mathrm{\;\!v}^{t\,\mathrm{mod}}\!-(1\!-\!B)\!\left[\tenq{\varepsilon}_\mathrm{\;\!p}^{t}+\Delta\tilde{\bar{\lambda}}^{t+1}\,\ProjectionTensDev:\tenq{\tilde{\sigma}}^{t+1}\right]\right]\,,
\end{equation}
which is an implicit, however, linear relation for the stress at the end of the increment $\tenq{\tilde{\sigma}}^{t+1}$ and can therefore be cast into explicit form, without solving any system of equations
\begin{equation}\label{stress_alg_2_exp}
		\tenq{\tilde{\sigma}}^{t+1}\!=\!(1-\tilde{d}^{\;\!t+1})\!\left[3\,\kappa\,\ProjectionTensSph\!+\!\dfrac{2\mu\ \ProjectionTensDev}{1\!+\![1-\tilde{d}^{\;\!t+1}](1\!-\!B)\,2\mu\,\Delta\tilde{\bar{\lambda}}^{t+1}}\right]\!\!:\!\!\left[\!(1\!-\!B)\tenq{\varepsilon}^{t+1}\!-\!(1\!-\!B)\tenq{\varepsilon}_\mathrm{\;\!p}^{t}\!-\!\tenq{\varepsilon}_\mathrm{\;\!v}^{t\,\mathrm{mod}}\right].
\end{equation}
The necessary algorithmic tangent for the global N\textsc{ewton} algorithm becomes
\begin{equation}
	\tenq[2]{{C}}\phantom{}_\mathrm{t}=\dfrac{\mathrm{d\,}\tenq{\tilde{\sigma}}^{t+1}}{\mathrm{d\,}\tenq{\varepsilon}^{t+1}}=(1-\tilde{d}^{\;\!t+1})(1-B)\left[3\,\kappa\,\ProjectionTensSph+\dfrac{2\mu\ \ProjectionTensDev}{1+(1-\tilde{d}^{\;\!t+1})(1-B)\,2\mu\,\Delta\tilde{\bar{\lambda}}^{t+1}}\right]\ .
\end{equation}
From equation \eqref{stress_alg_2_exp} it can be observed, that the stress $\tenq{\tilde{\sigma}}^{t+1}$ is a linear function of the strain $\tenq{\varepsilon}^{t+1}$, therefore, the algorithmic tangent $\tenq[2]{{C}}_\mathrm{t}^{\;\!t+1}$\vspace{0.3ex} is timestep-wise constant and the global problem is---in absence of further sources of nonlinearity---linear, hence, the problem can be solved in a single step. 
\subsection{Anisotropic damage and inelasticity in unidirectional composites}\label{chap:yarn_material}
\subsubsection*{Continuous formulation}
In woven composites, a bundled assembly of fibers is referred to as a yarn \cite{Praud_2018_dis}. Between these unidirectional oriented fibers the matrix material is situated. Since one fiber bundle consists of possibly hundreds of these fibers, 
\begin{figure}[!h]
	\centering
	\includegraphics[width=0.77\textwidth]{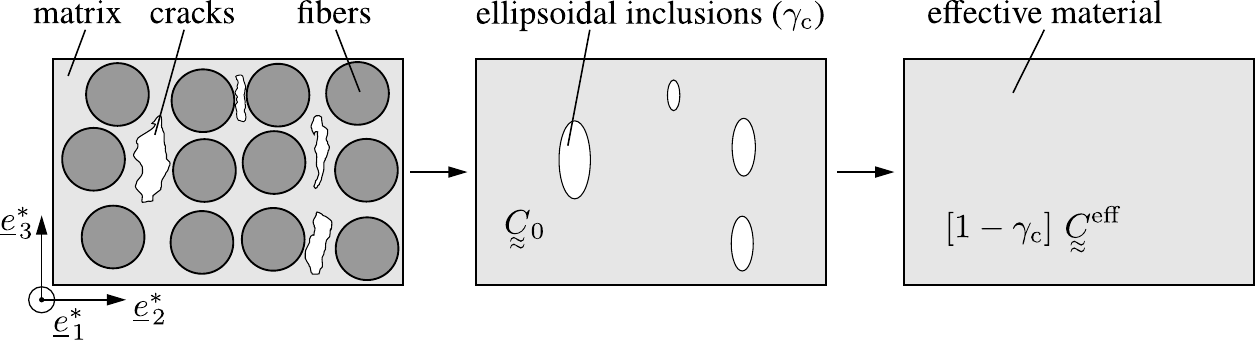}
	\caption{left: Yarn bundle with micro cracks, middle: micromechanical model with the cracks approximated by ellipsoidal inclusions, right: effective material.}
	\label{fig:Yarnscheme}
\end{figure}
it is favorable to find a homogenized constitutive model for the overall behavior of the yarns. The main modeling challenges are to describe the anisotropic stiffness, the evolution of micro cracks and the debonding of the fibers from the matrix, the rate-dependent behavior stemming from the matrix and the fiber failure among others. \citet{praud_2017_yarns} developed a model that incorporates anisotropic continuum damage and inelasticity induced by micro cracks, that is capable of describing at least the first two mentioned effects in a mean sense, while it is micromechanically motivated and enables a straightforward way of calibrating the model parameters with only a few experiments. The theory shall be explained only shortly before the numerical treatment is presented.

The model, as depicted in Figure\ref{fig:Yarnscheme}, assumes the existence of micro cracks with volume fraction $\gamma_\mathrm{c}\ge0$ and a priori known shape, e.g., through experiments, inside the virgin part of the material with a homogeneous anisotropic elastic behavior described with the elasticity tensor $\tenq[2]{{C}}\phantom{}_0$. The right-handed, orthonormal basis system $\mathcal{BS}=\left\{\underline{e}_{\;\!1}^*,\underline{e}_{\;\!2}^*,\underline{e}_{\;\!3}^*\right\}$ is chosen such that $\underline{e}_{\;\!1}^*$ is the fiber direction and $\underline{e}_{\;\!2}^*$ lies in the plane of the thin-walled woven composite structure. For a two-phase microstructure, the relation between the average phase strains can be expressed with localization tensors $\tenq[2]{{A}}\phantom{}_r$ \cite{chatzigeorgiou_2022}, which link an average phase strain $\tenq{\varepsilon}_{\;\!r}$ and the total strain through $\tenq{\varepsilon}_{\;\!r}=\tenq[2]{{A}}\phantom{}_r:\tenq{\varepsilon}$. The effective yarn strain is a sum of the average strains in the virgin part of the composite $\tenq{\varepsilon}_{\;\!0}$ and a fictitious crack strain $\tenq{\varepsilon}_\mathrm{\;\!c}$
\begin{equation}
	\tenq{\varepsilon}=\left(1-\gamma_\mathrm{c}\right)\underbrace{\tenq[2]{{A}}\phantom{}_0(\gamma_\mathrm{c}):\tenq{\varepsilon}}_{\tenq{\varepsilon}_{\;\!0}}+\gamma_\mathrm{c}\underbrace{\tenq[2]{{A}}\phantom{}_\mathrm{c}(\gamma_\mathrm{c}):\tenq{\varepsilon}}_{\tenq{\varepsilon}_\mathrm{\;\!c}}\ .
\end{equation}
The stress $\tenq{\sigma}_\mathrm{c}$ and strain $\tenq{\varepsilon}_\mathrm{\;\!c}$ inside the cracks have obviously no meaning, because there is no material and the overall yarn stress simply becomes
\begin{equation}
	\tenq{\sigma}=\left(1-\gamma_\mathrm{c}\right)\tenq{\sigma}_0\ .
\end{equation} 
The here applied \textsc{Mori}-\textsc{Tanaka} method \cite{Mori_1973} assumes, that far enough from the defects, the strain field can be replaced by its \enquote{matrix} (here virgin material) volume average, leading to the strain concentration tensors
\begin{equation}
	\tenq[2]{{A}}\phantom{}_0(\gamma_\mathrm{c})=\left[\tenq[2]{{I}}+\gamma_\mathrm{c}\,\left[\tenq[2]{{T}}\phantom{}_\mathrm{c}-\tenq[2]{{I}}\right]\right]^{-1}\ \ \ \text{and}\ \ \ \ \tenq[2]{{A}}\phantom{}_\mathrm{c}(\gamma_\mathrm{c})=\tenq[2]{{T}}\phantom{}_\mathrm{c}:\tenq[2]{{A}}\phantom{}_0(\gamma_\mathrm{c})\ .
\end{equation}
Therein, the interaction tensor $\tenq[2]{{T}}\phantom{}_\mathrm{c}$, which can be computed with aid of the \textsc{Eshelby} tensor~$\tenq[2]{{S}}_\mathrm{E}$~\cite{Eshelby_1957} that has to be evaluated numerically in the here relevant case of an anisotropic surrounding media \cite{Gavazzi_1990}, reads
\begin{equation}
	\tenq[2]{{T}}\phantom{}_\mathrm{c}=\left[\tenq[2]{{I}}-\tenq[2]{{S}}\phantom{}_\mathrm{E}\right]^{-1}\ .
\end{equation}
The total strain $\tenq{\varepsilon}$ is assumed to be additively split into an elastic strain $\tenq{\varepsilon}_\mathrm{\;\!e}$ and an inelastic strain $\tenq{\varepsilon}_\mathrm{\;\!s}$, that occurs as a result of the damage mechanisms and remains permanently as a consequence of the nonclosure of the microcracks
\begin{equation}
	\tenq{\varepsilon}=\tenq{\varepsilon}_\mathrm{\;\!e}+\tenq{\varepsilon}_\mathrm{\;\!s}\ \ \ \rightarrow\ \ \ \tenq{\varepsilon}_\mathrm{\;\!e}=\tenq{\varepsilon}-\tenq{\varepsilon}_\mathrm{\;\!s}\ .
\end{equation}
Therewith, the overall stress can be computed from the stress in the virgin part of the material, giving
\begin{gather}\label{stress_state_yarn}
	\tenq{\sigma}
    =\!\left(1-\gamma_\mathrm{c}\right)\tenq[2]{{C}}\phantom{}_0\!:\!\tenq{\varepsilon}_{\mathrm{\;\!e}0}
    \!=\!\left(1-\gamma_\mathrm{c}\right)\underbrace{\tenq[2]{{C}}\phantom{}_0\!:\!\tenq[2]{{A}}\phantom{}_0(\gamma_\mathrm{c})}_{:=\tenq[2]{{C}}^\mathrm{eff}(\gamma_\mathrm{c})}\!:\!\tenq{\varepsilon}_\mathrm{\;\!e}
    \!=\!\left(1-\gamma_\mathrm{c}\right)\tenq[2]{{C}}^\mathrm{eff}(\gamma_\mathrm{c})\!:\!\left[\;\!\tenq{\varepsilon}-\tenq{\varepsilon}_\mathrm{\;\!s}\;\!\right]\,.\vspace{-1.5ex}\raisetag{3ex}
\end{gather}
The evolution of the micro crack volume fraction $\gamma_\mathrm{c}$ is chosen in the form of a \textsc{Weibull}-like relation, with parameters $\gamma_\mathrm{c}^\infty$, $S$ and $\beta$, as a function of a crack evolution activation $H_\mathrm{c}$
\begin{equation}\label{gamma_c_evo}
	\gamma_\mathrm{c}=\gamma_\mathrm{c}^\infty\left[1-\mathrm{exp}\left(-\left[\dfrac{\langle\mathrm{sup}(H_\mathrm{c})-1\rangle}{S}\right]^\beta\right)\right]\ .
\end{equation}
The cracks will only evolve if the criterion $H_\mathrm{c}$ takes values above its ever encountered value, and consequently the crack density cannot decline. The function $H_\mathrm{c}$ 
\begin{gather}
	H_\mathrm{c}=\sqrt{\tenq{\sigma}_0:\tenq[2]{{H}}:\tenq{\sigma}_0}=\dfrac{1}{1-\gamma_\mathrm{c}}\sqrt{\tenq{\sigma}:\tenq[2]{{H}}:\tenq{\sigma}}\ ,\ \ \left[\tenq[2]{{H}}\,\right]=
	\begin{bmatrix}
		 0 & 0 & 0 & 0 & 0 & 0 \\
		  & 1/R_{22}^2 & 0 & 0 & 0 & 0 \\
		  &  & 0 & 0 & 0 & 0 \\
		  &  &  & 0 & 0 & 0 \\
		  &  &  &  & 0 & 0 \\
		 \mathrm{sym} &  &  &  &  & 1/R_{12}^2
	\end{bmatrix}_{\mathcal{BS}}\raisetag{-1ex}
\end{gather}
is formulated in the stresses and set up such, that it is only sensitive to in plane shear and traverse tension. The evolution of the permanent inelastic strains is coupled to the evolution of the micro crack volume fraction with aid of another function $H_\mathrm{s}$ by
\begin{gather}\label{anel_evo}
	\tenq{\dot{\varepsilon}}_\mathrm{\;\!s}=\dfrac{\tenq[2]{{F}}:\tenq{\sigma}}{H_s}\dot{\gamma}_\mathrm{c}\ \ \text{with}\ \ H_\mathrm{s}=\sqrt{\tenq{\sigma}:\tenq[2]{{F}}:\tenq{\sigma}}\ \ \text{and}\ \ \left[\tenq[2]{{F}}\,\right]=
	\begin{bmatrix}
		0 & 0 & 0 & 0 & 0 & 0 \\
		& a_{22}^2 & 0 & 0 & 0 & 0 \\
		&  & 0 & 0 & 0 & 0 \\
		&  &  & 0 & 0 & 0 \\
		&  &  &  & 0 & 0 \\
		\mathrm{sym} &  &  &  &  & a_{12}^2
	\end{bmatrix}_{\mathcal{BS}}\ .\raisetag{-1ex}
\end{gather}
The criterion $H_\mathrm{s}$ is again anisotropic in the stress and  only the inelastic strains ${\varepsilon}_{\mathrm{s}\,22}$ and ${\varepsilon}_{\mathrm{s}\,12}(={\varepsilon}_{\mathrm{s}\,21}$) will evolve.
\subsubsection*{Implicit integration}
In the works of \citet{praud_2017_yarns} an implicit cutting-plane algorithm was deployed to solve the nonlinear equations. This is here however unfavorable, because the viscoelastic-viscplastic-damage matrix law as described in Section \ref{chap:matrix_material} was already discretized using an ImplEx method, yielding a globally timestep-wise linear problem. This makes it favorable to use an ImplEx algorithm also for the yarn material formulation.

The implicitly integrated evolution equations are solved using a \textsc{Newton} method. The evolution equations \eqref{gamma_c_evo} and \eqref{anel_evo} read in discretized form
\begin{equation}\label{gamma_c_evo_dis}
	\gamma_\mathrm{c}^{t+1}=\gamma_\mathrm{c}^\infty\Big[1-\mathrm{exp}(-\left[\dfrac{\langle\mathrm{sup}(H_\mathrm{c}(\tenq{\sigma}^{t+1}))-1\rangle}{S}\right]^\beta)\Big]
\end{equation}
and
\begin{equation}\label{anel_evo_dis}
	\tenq{\varepsilon}_\mathrm{\;\!s}^{t+1}=\tenq{\varepsilon}_\mathrm{\;\!s}^{t}+\dfrac{\tenq[2]{{F}}:\tenq{\sigma}^{t+1}}{H_s(\tenq{\sigma}^{t+1})}\left[\gamma_\mathrm{c}^{t+1}-\gamma_\mathrm{c}^{t}\right]
\end{equation}
The overall stress in time discretized form becomes with equation \eqref{stress_state_yarn}
\begin{equation}\label{stress_state_yarn_disc}
	\tenq{\sigma}^{t+1}=\left[1-\gamma_\mathrm{c}^{t+1}\right]\tenq[2]{{C}}^\mathrm{eff}(\gamma_\mathrm{c}^{t+1}):\left[\tenq{\varepsilon}^{t+1}-\tenq{\varepsilon}_\mathrm{\;\!s}^{t+1}\right]\ .
\end{equation}
In a trial step it has to be assumed that no crack grows occur, i.e., $\gamma_\mathrm{c}^{t+1}=\gamma_\mathrm{c}^{t}$ and $\tenq{\varepsilon}_\mathrm{\;\!s}^{t+1}=\tenq{\varepsilon}_\mathrm{\;\!s}^{t}$. Then the activation criterion $H_\mathrm{c}(\tenq{\sigma}^{t+1})$ has to be evaluated. If $H_\mathrm{c}(\tenq{\sigma}^{t+1})$ is greater than one and $\mathrm{sup}(H_\mathrm{c})$, a correction must be performed, otherwise the trial step was correct. In order to solve the equations efficiently, equation \eqref{anel_evo_dis} has to be inserted into \eqref{stress_state_yarn_disc}
\begin{equation}\label{stress_state_yarn_disc_mod}
	\tenq{\sigma}^{t+1}=\left[1-\gamma_\mathrm{c}^{t+1}\right]\tenq[2]{{C}}^\mathrm{eff}(\gamma_\mathrm{c}^{t+1}):\left[\tenq{\varepsilon}^{t+1}-\tenq{\varepsilon}_\mathrm{\;\!s}^{t}-\dfrac{\tenq[2]{{F}}:\tenq{\sigma}^{t+1}}{H_s(\tenq{\sigma}^{t+1})}\left[\gamma_\mathrm{c}^{t+1}-\gamma_\mathrm{c}^{t}\right]\right]\ .
\end{equation}
It can be observed that only $\gamma_\mathrm{c}^{t+1}$, $\sigma_2^{t+1}$ and $\sigma_6^{t+1}$ are coupled through equations \eqref{stress_state_yarn_disc_mod} and \eqref{gamma_c_evo_dis}. Switching to V\textsc{oigt} notation, the following residual $\underline{R}$ is set up from equations \eqref{stress_state_yarn_disc_mod} and \eqref{gamma_c_evo_dis}.
\begin{equation}
	\underline{R}:=
	\begin{bmatrix}
		\vspace{1ex}
		R_{\sigma_2}\\[1.0ex]
		R_{\sigma_6}\\[3.0ex]
		R_{\gamma_\mathrm{c}}
		\vspace{1ex}
	\end{bmatrix}
	=
	\begin{bmatrix}
		\sigma_2^{t+1}-\left[1-\gamma_\mathrm{c}^{t+1}\right]\left[\tenq[2]{C}^\mathrm{eff}\right]_{2:}^\text{\footnotemark[2]}(\gamma_\mathrm{c}^{t+1})\cdot\left[\left[\tenq{\varepsilon}^{t+1}\right]-\left[\tenq{\varepsilon}_\mathrm{\;\!s}^{t+1}\right](\gamma_\mathrm{c}^{t+1},\sigma_2^{t+1},\sigma_6^{t+1})\right]\\[1.5ex]
		\sigma_6^{t+1}-\left[1-\gamma_\mathrm{c}^{t+1}\right]\left[\tenq[2]{C}^\mathrm{eff}\right]_{6:}(\gamma_\mathrm{c}^{t+1})\cdot\left[\left[\tenq{\varepsilon}^{t+1}\right]-\left[\tenq{\varepsilon}_\mathrm{\;\!s}^{t+1}\right](\gamma_\mathrm{c}^{t+1},\sigma_2^{t+1},\sigma_6^{t+1})\right]\\[0.5ex]
		\gamma_\mathrm{c}^{t+1}-\gamma_\mathrm{c}^\infty\Big[1-\mathrm{exp}(-\left[\dfrac{ H_\mathrm{c}(\gamma_\mathrm{c}^{t+1},\sigma_2^{t+1},\sigma_6^{t+1})-1}{S}\right]^\beta)\Big]
	\end{bmatrix}
\end{equation}\footnotetext[2]{The writing $\left[\tenq[2]{C}^\mathrm{eff}\right]_{2:}$ means in this context the second matrix row of the V\textsc{oigt} representation of $\tenq[2]{C}^\mathrm{eff}$.}
The N\textsc{ewton} solution method reads for the problem
\begin{equation}
	\begin{bmatrix}
		K_{\sigma_2\,\sigma_2}&K_{\sigma_2\,\sigma_6}&K_{\sigma_2\,\gamma_\mathrm{c}}\\
		K_{\sigma_6\,\sigma_2}&K_{\sigma_6\,\sigma_6}&K_{\sigma_6\,\gamma_\mathrm{c}}\\
		K_{\gamma_\mathrm{c}\,\sigma_2}&K_{\gamma_\mathrm{c}\,\sigma_6}&K_{\gamma_\mathrm{c}\,\gamma_\mathrm{c}}
	\end{bmatrix}_k
	\cdot
	\begin{bmatrix}
		\Delta \sigma_2\\
		\Delta \sigma_6\\
		\Delta \gamma_\mathrm{c}
	\end{bmatrix}_{k+1}
	=
	-\underline{R}_{\,k}
	\ \ \ \ \text{and}\ \ \ \ \ 
	\begin{matrix}
		\sigma_2^{k+1}=\sigma_2^{k}+\Delta \sigma_2^{k+1}\\[0.5ex]
		\sigma_6^{k+1}=\sigma_6^{k}+\Delta \sigma_6^{k+1}\\[0.5ex]
		\gamma_\mathrm{c}^{k+1}=\gamma_\mathrm{c}^{k}+\Delta \gamma_\mathrm{c}^{k+1} \,.
	\end{matrix}
\end{equation}
Therein, the necessary derivatives for the N\textsc{ewton} scheme read
\begin{equation}
	K_{\sigma_2\,\sigma_2}^k=1+\left[1-\gamma_\mathrm{c}^k\right]\cdot \left[C^\mathrm{eff}_{22}\cdot\dfrac{\mathrm{d}\,\varepsilon_{\mathrm{s}\,2}}{\mathrm{d}\,\sigma_2^{k}}+C^\mathrm{eff}_{26}\cdot\dfrac{\mathrm{d}\,\varepsilon_{\mathrm{s}\,6}}{\mathrm{d}\,\sigma_2^{k}}\right]\ ,
\end{equation}
\begin{equation}
	K_{\sigma_2\,\sigma_6}^k=\left[1-\gamma_\mathrm{c}^k\right]\cdot \left[C^\mathrm{eff}_{22}\cdot\dfrac{\mathrm{d}\,\varepsilon_{\mathrm{s}\,2}}{\mathrm{d}\,\sigma_6^{k}}+C^\mathrm{eff}_{26}\cdot\dfrac{\mathrm{d}\,\varepsilon_{\mathrm{s}\,6}}{\mathrm{d}\,\sigma_6^{k}}\right]\ ,
\end{equation}
\begin{align}
	\nonumber K_{\sigma_2\,\gamma_\mathrm{c}}^k=&\left[\tenq[2]{C}^\mathrm{eff}\right]_{2:}\cdot\left[\left[\tenq{\varepsilon}^{t+1}\right]-\left[\tenq{\varepsilon}_{\;\!s}\right]\right]+\left[1-\gamma_\mathrm{c}^k\right]\cdot\left[C^\mathrm{eff}_{22}\cdot\dfrac{\mathrm{d}\,\varepsilon_{\mathrm{s}\,2}}{\mathrm{d}\,\gamma_\mathrm{c}^k}+C^\mathrm{eff}_{26}\cdot\dfrac{\mathrm{d}\,\varepsilon_{\mathrm{s}\,6}}{\mathrm{d}\,\gamma_\mathrm{c}^k}\right]\\
	&+\left[1-\gamma_\mathrm{c}^k\right]\ \left[\tenq[2]{C}^\mathrm{eff}\right]_{2:}\cdot\left[\tenq[2]{A}_{\,0}\right]\cdot\left[\left[\tenq[2]{T}_{\,\mathrm{c}}\right]-\left[\tenq[2]{I}\right]\right]\cdot\left[\left[\tenq{\varepsilon}^{t+1}\right]-\left[\tenq{\varepsilon}_{\;\!s}\right]\right]\ ,
\end{align}
\begin{equation}
	K_{\sigma_6\,\sigma_2}^k=\left[1-\gamma_\mathrm{c}^k\right]\cdot \left[C^\mathrm{eff}_{62}\cdot\dfrac{\mathrm{d}\,\varepsilon_{\mathrm{s}\,2}}{\mathrm{d}\,\sigma_2^{k}}+C^\mathrm{eff}_{66}\cdot\dfrac{\mathrm{d}\,\varepsilon_{\mathrm{s}\,6}}{\mathrm{d}\,\sigma_2^{k}}\right]\ ,
\end{equation}
\begin{equation}
	K_{\sigma_6\,\sigma_6}^k=1+\left[1-\gamma_\mathrm{c}^k\right]\cdot \left[C^\mathrm{eff}_{62}\cdot\dfrac{\mathrm{d}\,\varepsilon_{\mathrm{s}\,2}}{\mathrm{d}\,\sigma_6^{k}}+C^\mathrm{eff}_{66}\cdot\dfrac{\mathrm{d}\,\varepsilon_{\mathrm{s}\,6}}{\mathrm{d}\,\sigma_6^{k}}\right]\ ,
\end{equation}
\begin{align}
	\nonumber K_{\sigma_6\,\gamma_\mathrm{c}}^k=&\left[\tenq[2]{C}^\mathrm{eff}\right]_{6:}\cdot\left[\left[\tenq{\varepsilon}^{t+1}\right]-\left[\tenq{\varepsilon}_{\;\!s}\right]\right]+\left[1-\gamma_\mathrm{c}^k\right]\cdot\left[C^\mathrm{eff}_{62}\cdot\dfrac{\mathrm{d}\,\varepsilon_{\mathrm{s}\,2}}{\mathrm{d}\,\gamma_\mathrm{c}^k}+C^\mathrm{eff}_{66}\cdot\dfrac{\mathrm{d}\,\varepsilon_{\mathrm{s}\,6}}{\mathrm{d}\,\gamma_\mathrm{c}^k}\right]\\
	&+\left[1-\gamma_\mathrm{c}^k\right]\cdot\left[\tenq[2]{C}^\mathrm{eff}\right]_{6:}\cdot\left[\tenq[2]{A}_{\,0}\right]\cdot\left[\left[\tenq[2]{T}_{\,\mathrm{c}}\right]-\left[\tenq[2]{I}\right]\right]\cdot\left[\left[\tenq{\varepsilon}^{t+1}\right]-\left[\tenq{\varepsilon}_{\;\!s}\right]\right]\ ,
\end{align}
\begin{equation}
	K_{\gamma_\mathrm{c}\,\sigma_2}^k=-\gamma_\mathrm{c}^\infty\exp(-\left[\dfrac{H_\mathrm{c}-1}{S}\right]^\beta)\cdot\left[\dfrac{H_\mathrm{c}-1}{S}\right]^{2\beta-1}\cdot\dfrac{1}{S}\cdot\dfrac{\mathrm{d}\,H_\mathrm{c}}{\mathrm{d}\,\sigma_2^{k}}\ ,
\end{equation}
\begin{equation}
	K_{\gamma_\mathrm{c}\,\sigma_6}^k=-\gamma_\mathrm{c}^\infty\exp(-\left[\dfrac{H_\mathrm{c}-1}{S}\right]^\beta)\cdot\left[\dfrac{H_\mathrm{c}-1}{S}\right]^{2\beta-1}\cdot\dfrac{1}{S}\cdot\dfrac{\mathrm{d}\,H_\mathrm{c}}{\mathrm{d}\,\sigma_6^{k}}\ ,
\end{equation}
\begin{equation}
	K_{\gamma_\mathrm{c}\,\gamma_\mathrm{c}}^k=1-\gamma_\mathrm{c}^\infty\exp(-\left[\dfrac{H_\mathrm{c}-1}{S}\right]^\beta)\cdot\left[\dfrac{H_\mathrm{c}-1}{S}\right]^{2\beta-1}\cdot\dfrac{1}{S}\cdot\dfrac{\mathrm{d}\,H_\mathrm{c}}{\mathrm{d}\,\gamma_\mathrm{c}^{k}}\ ,
\end{equation}
\begin{equation}
	\dfrac{\mathrm{d}\,\varepsilon_{\mathrm{s}\,2}}{\mathrm{d}\,\sigma_2^{k}}=\left[\gamma_\mathrm{c}^k-\gamma_\mathrm{c}^t\right]a_{22}^2\,\dfrac{H_\mathrm{s}-\sigma_2^{k}\cdot\dfrac{\mathrm{d}\,H_\mathrm{s}}{\mathrm{d}\,\sigma_2^{k}}}{H_\mathrm{s}^2}\ ,\ \ \ \dfrac{\mathrm{d}\,\varepsilon_{\mathrm{s}\,6}}{\mathrm{d}\,\sigma_2^{k}}=-\left[\gamma_\mathrm{c}^k-\gamma_\mathrm{c}^t\right]a_{12}^2\,\sigma_6^{k}\dfrac{\dfrac{\mathrm{d}\,H_\mathrm{s}}{\mathrm{d}\,\sigma_2^{k}}}{H_\mathrm{s}^2}\ ,\ \ \dfrac{\mathrm{d}\,H_\mathrm{s}}{\mathrm{d}\,\sigma_2^{k}}=\dfrac{a_{22}^2\,\sigma_2^{k}}{H_\mathrm{s}}\ ,
\end{equation}
\begin{equation}
	\dfrac{\mathrm{d}\,\varepsilon_{\mathrm{s}\,2}}{\mathrm{d}\,\sigma_6^{k}}=-\left[\gamma_\mathrm{c}^k-\gamma_\mathrm{c}^t\right]a_{22}^2\,\sigma_2^{k}\,\dfrac{\dfrac{\mathrm{d}\,H_\mathrm{s}}{\mathrm{d}\,\sigma_6^{k}}}{H_\mathrm{s}^2}\ ,\ \ \ \dfrac{\mathrm{d}\,\varepsilon_{\mathrm{s}\,6}}{\mathrm{d}\,\sigma_6^{k}}=\left[\gamma_\mathrm{c}^k-\gamma_\mathrm{c}^t\right]a_{12}^2\,\dfrac{H_\mathrm{s}-\sigma_6^{k}\dfrac{\mathrm{d}\,H_\mathrm{s}}{\mathrm{d}\,\sigma_6^{k}}}{H_\mathrm{s}^2}\ ,\ \ \dfrac{\mathrm{d}\,H_\mathrm{s}}{\mathrm{d}\,\sigma_6^{k}}=\dfrac{a_{12}^2\,\sigma_6^{k}}{H_\mathrm{s}}\ ,
\end{equation}
\begin{equation}
	\dfrac{\mathrm{d}\,\varepsilon_{\mathrm{s}\,2}}{\mathrm{d}\,\gamma_\mathrm{c}^k}=\dfrac{\sigma^k_2a_{22}^2}{H_\mathrm{s}}\ ,\ \ \ \dfrac{\mathrm{d}\,\varepsilon_{s\,6}}{\mathrm{d}\,\gamma_\mathrm{c}^k}=\dfrac{\sigma^k_6a_{12}^2}{H_\mathrm{s}}\ ,
\end{equation}
and 
\begin{equation}
	\dfrac{\mathrm{d}\,H_\mathrm{c}}{\mathrm{d}\,\sigma_2^{k}}=\dfrac{1}{\left[1-\gamma_\mathrm{c}^k\right]^2\cdot H_\mathrm{c}}\cdot\dfrac{\sigma^k_2}{R_{22}^2}\ ,\ \ \dfrac{\mathrm{d}\,H_\mathrm{c}}{\mathrm{d}\,\sigma_6^{k}}=\dfrac{1}{\left[1-\gamma_\mathrm{c}^k\right]^2\cdot H_\mathrm{c}}\cdot\dfrac{\sigma^k_6}{R_{12}^2}\ ,\ \ \dfrac{\mathrm{d}\,H_\mathrm{c}}{\mathrm{d}\,\gamma_\mathrm{c}^{k}}=H_\mathrm{c}\cdot\dfrac{1}{1-\gamma_\mathrm{c}^k}\ .
\end{equation}
The criterion for convergence is defined by
\begin{equation}
	\dfrac{|R_{\sigma_2}|+|R_{\sigma_6}|}{\max(\left|\underline{\sigma}\right|)}+\dfrac{|R_{\gamma_\mathrm{c}}|}{\gamma_\mathrm{c}}<\varepsilon\ .
\end{equation}

\subsubsection*{Explicit integration for ImplEx}
The explicit measures of the material formulation, whose stress and algorithmic tangent are to be returned to the element routine, are marked with a tilde $(\tilde{\phantom{a}})$. The measures from the last step without tilde are those computed with the \textsc{Euler} \emph{backward} method. In order to compute the stress at the end of the increment through equation \eqref{stress_state_yarn}, it is necessary to know $\tilde{\gamma}_\mathrm{c}^{t+1}$ and  $\tenq{\tilde{\varepsilon}}_{\;\!s}^{\,t+1}$. The crack volume fraction is explicitly integrated by
\begin{equation}\label{gamma_exp}
	\tilde{\gamma}_\mathrm{c}^{t+1}=\gamma^{t}_\mathrm{c}+\Delta t\,\dot{\gamma}_\mathrm{c}^{t}=\gamma_\mathrm{c}^{t}+\Delta t\,\dfrac{\gamma_\mathrm{c}^{t}-\gamma_\mathrm{c}^{t-1}}{t^t-t^{t-1}}\ .
\end{equation}
The inelastic strains are integrated analogously to equation \eqref{anel_evo_dis} as
\begin{equation}\label{anel_exp}
	\tenq{\tilde{\varepsilon}}_\mathrm{\;\!s}^{\,t+1}=\tenq{\varepsilon}_\mathrm{\;\!s}^{t}+\dfrac{\tenq[2]{{F}}:\tenq{\tilde{\sigma}}^{t+1}}{\tilde{H}_\mathrm{s}^{t+1}}\left[\tilde{\gamma}_\mathrm{c}^{t+1}-\gamma_\mathrm{c}^{t}\right]\ \ \text{with}\ \ \ \tilde{H}_\mathrm{s}^{t+1}=H_\mathrm{s}^{t}+\Delta t\,\dot{H}_\mathrm{s}^{t}=H_\mathrm{s}^{t}+\dfrac{H_\mathrm{s}^{t}-H_\mathrm{s}^{t-1}}{t^t-t^{t-1}}
\end{equation}
wherein $\tenq{\tilde{\sigma}}^{t+1}$ is inserted as the stress measure, but $\tilde{H}_\mathrm{s}^{t+1}$ is integrated explicitly, since it is a nonlinear function of the stress and therefore an iterative evaluation can be avoided. The equations \eqref{stress_state_yarn}, \eqref{gamma_exp} and \eqref{anel_exp} together yield
\begin{equation}
	\tenq{\tilde{\sigma}}^{t+1}=\left(1-\tilde{\gamma}_\mathrm{c}^{t+1}\right)\tenq[2]{{C}}^\mathrm{eff}(\tilde{\gamma}_\mathrm{c}^{t+1}):\left[\tenq{\varepsilon}^{t+1}-\tenq{\varepsilon}_\mathrm{\;\!s}^{t}+\dfrac{\tenq[2]{{F}}:\tenq{\tilde{\sigma}}^{t+1}}{\tilde{H}_\mathrm{s}^{t+1}}\left[\tilde{\gamma}_\mathrm{c}^{t+1}-\gamma_\mathrm{c}^{t}\right]\right]\ ,
\end{equation}
which is an implicit equation for the stress $\tenq{\tilde{\sigma}}^{t+1}$, but can be written explicitly as
\begin{equation}\label{stress_yarn_final}
	\tenq{\tilde{\sigma}}^{t+1}\!=\!\left(1-\tilde{\gamma}_\mathrm{c}^{t+1}\right)\left[\tenq[2]{{I}}-\dfrac{\left[\tilde{\gamma}_\mathrm{c}^{t+1}-\gamma_\mathrm{c}^{t}\right]\left[1-\tilde{\gamma}_\mathrm{c}^{t+1}\right]}{\tilde{H}_\mathrm{s}^{t+1}}\ \tenq[2]{{C}}^\mathrm{eff}(\tilde{\gamma}_\mathrm{c}^{t+1}):\tenq[2]{{F}}\ \right]^{-1}:\ \tenq[2]{{C}}^\mathrm{eff}(\tilde{\gamma}_\mathrm{c}^{t+1}):\left[\tenq{\varepsilon}^{t+1}-\tenq{\varepsilon}_\mathrm{\;\!s}^{t}\right]
\end{equation}
and gives the algorithmic consistent tangent
\begin{equation}
	\tenq[2]{{C}}\phantom{}_\mathrm{t}=\dfrac{\mathrm{d\,}\tenq{\tilde{\sigma}}^{t+1}}{\mathrm{d\,}\tenq{\varepsilon}^{t+1}}=\left(1-\tilde{\gamma}_\mathrm{c}^{t+1}\right)\left[\tenq[2]{{I}}-\dfrac{\left[\tilde{\gamma}_\mathrm{c}^{t+1}-\gamma_\mathrm{c}^{t}\right]\left[1-\tilde{\gamma}_\mathrm{c}^{t+1}\right]}{\tilde{H}_\mathrm{s}^{t+1}}\ \tenq[2]{{C}}^\mathrm{eff}(\tilde{\gamma}_\mathrm{c}^{t+1}):\tenq[2]{{F}}\ \right]^{-1}:\ \tenq[2]{{C}}^\mathrm{eff}(\tilde{\gamma}_\mathrm{c}^{t+1})\ .
\end{equation}
Since the explicit stress $\tenq{\tilde{\sigma}}^{t+1}$ is a linear function of the strain over one time step as can be seen from equation \eqref{stress_yarn_final}, the global \textsc{Newton} scheme will---in absence from other sources of nonlinearity---converge in one increment.

\end{document}